\RequirePackage{fix-cm}
\documentclass[oneside,a4paper,11pt,numbers=noenddot,headinclude,footinclude=true,cleardoublepage=empty]{scrartcl}

\usepackage[main=english]{babel}
\usepackage[T1]{fontenc}

\usepackage[labelfont=bf]{caption}

\usepackage{graphicx}
\usepackage{subcaption}
\usepackage{mwe}

\usepackage{amsmath,amsfonts,amssymb,amsthm,booktabs,color,epsfig,graphicx,hyperref,url,mathtools}

\theoremstyle{plain}% Theorem-like structures provided by amsthm.sty
\newtheorem{theorem}{Theorem}

\newtheorem{proposition}{Proposition}
\newtheorem{lemma}{Lemma}
\newtheorem{corollary}{Corollary}

\theoremstyle{definition}

\newtheorem{remark}{Remark}

\DeclareMathOperator*{\argmin}{arg\,min}

\usepackage{geometry}
\usepackage{array}
\newcolumntype{L}[1]{>{\raggedright\let\newline\\\arraybackslash\hspace{0pt}}m{#1}}
\newcolumntype{C}[1]{>{\centering\let\newline\\\arraybackslash\hspace{0pt}}m{#1}}
\newcolumntype{R}[1]{>{\raggedleft\let\newline\\\arraybackslash\hspace{0pt}}m{#1}}

\usepackage[dvipsnames]{xcolor}

\title{Robust signal recovery in hadamard spaces}
\author{Georg Köstenberger, Thomas Stark}

\begin{document}
\frenchspacing
\raggedbottom

\maketitle

%\author[1]{Georg Köstenberger\corref{mycorrespondingauthor}}
%\author[1]{Thomas Stark}

%\address[1]{Department of Statistics and Operations Research, University of Vienna, Oskar-Morgenstern-Platz 1, A-1090 Vienna}
%\address[2]{Department of Statistics and Operations Research, University of Vienna, Oskar-Morgenstern-Platz 1, A-1090 Vienna}

%\cortext[mycorrespondingauthor]{Corresponding author. Email address: \url{georg.koestenberger@univie.ac.at}}

\begin{abstract}
%Complete metric spaces of non-positive curvature (also known as Hadamard spaces) have seen considerable success in applications and theory alike. 
 % A basic problem in probability theory and statistics is the computation of averages or means.
  %There are, however, various non-equivalent notions of \textit{mean} in Hadamard spaces. 
  %Some are well-behaved in theory, others are computationally tractable, and to the best of our knowledge, none are both.
  %In many cases one is interested in the asymptotic behavior of such means.
  We analyze the stability of (strong) laws of large numbers in Hadamard spaces with respect to distributional perturbations.
  For the inductive means of a sequence of independent, but not necessarily identically distributed random variables, we provide a concentration inequality in quadratic mean, as well as a strong law of large numbers, generalizing a classical result of K.-T. Sturm. 
  For the Fr\'echet mean, we generalize H. Ziezold's law of large numbers in Hadamard spaces.
  In this case, we neither require our data to be independent, nor identically distributed; reasonably mild conditions on the first two moments of our sample are enough.
  Additionally, we look at data contamination via a model inspired by Huber's $\varepsilon$-contamination model, in which we replace a random portion of the data with noise.
  In the most general setup, we do neither require the data, nor the noise to be i.i.d., nor do we require the noise to be independent of the data.
  To analyze the stability of the (non-symmetric) inductive mean with respect to data loss, data permutation, and noise, a resampling scheme is introduced, and sufficient conditions for its convergence are provided.
  These results suggest that means in Hadamard spaces are as robust as in Euclidean spaces.
  This is underlined by a small simulation study, in which we compare the robustness of means on the manifold of positive definite matrices, with means on open books.

%---------------------------------

  %Classical statistical theory has been developed under the assumption that the data belongs to a linear space. 
  %However, in many applications the intrinsic geometry of the data is more intricate. 
  %Neglecting this frequently yields suboptimal or outright unusable results, e.g., taking the pixel-wise average \revrem{of images}{R1-NI-1}{of portraits} typically \revrem{results in noise}{R1-NI}{does not yield an image of a person}.
  %Incorporating the intrinsic geometry of a dataset into statistical analysis is a highly non-trivial task.
  %In fact different underlying geometries necessitate different approaches, and allow for results of varying strength.
  %Perhaps the most common non-linear geometries appearing in statistical applications are metric spaces of non-positive curvature, such as the manifold of symmetric, positive (semi-)definite matrices. 
  %In this paper we introduce a (strong) law of large numbers for independent, but not necessarily identically distributed random variables taking values in complete spaces of non-positive curvature. 
  %Using this law of large numbers, we justify a stochastic approximation scheme for the limit of Fr\'{e}chet means on such spaces.
  %Apart from rendering the problem of computing Fr\'{e}chet means computationally more tractable, the structure of this scheme suggests, that averaging operations on Hadamard spaces are more stable than previous results would suggest.
  \end{abstract}

\section{Introduction}
%With a general increase in computing power in recent years and the advent of big data, many statistical procedures previously banned to the realm of theory, have become computationally feasible.
%In particular, the resolution of many datasets has become so fine-grained, that one may think of it as a very high dimensional vector or function.
%Hence one is interested in a statistical theory for samples drawn from infinite dimensional data. 
%The most general setup under which such a theory is thinkable are metric spaces.

%Contemporary machine learning techniques are inherently non-linear, we routinely work with images, shapes, or general point clouds, and even if our data is sampled from a linear (typically high-dimensional) space, it is frequently assumed that the data points lie on a latent, lower-dimensional manifold, (see \cite{fefferman-manifold} for empirical verification and the history of this hypothesis) to overcome the curse of dimensionality.  
%Ignoring the non-linear, latent structure of the data can have fatal consequences for our data analysis, e.g. taking the pixel-wise average in a dataset of portraits, will typically not result in an average face.
In modern statistics, non-linear spaces play an increasingly important role, and the need for a non-linear statistical theory has been recognized in the literature long ago \cite{ziezold1977, slln-centroid-cpt-metric-space}.
One of the most general setups under which a reasonably complete non-linear statistical theory is feasible \cite{es-sahib-heinich,sturm,sturm-martingale,le-gouic-fast-convergence,brunel,escande2023}, are metric spaces of non-positive curvature in the sense of Alexandrov \cite{alexandrov57} (also known as CAT($0$) spaces).
These spaces include prominent examples, such as the manifold of positive definite, symmetric matrices, Hilbert spaces, metric trees, and Riemannian manifolds of non-positive sectional curvature, and appear in fields as diverse as computational biology \cite{Billera-phylogenetic-trees}, image processing and medicine \cite{pennec2019riemannian, med1, med2}, dynamical systems \cite{burago-ferleger-kononenko}, probability theory \cite{sturm}, statistics \cite{huber64, huber65,huckemann}, group theory \cite{gromov1987,qing-rafi-random-walks, gromov-random-walks} and, of course geometry, where CAT$(0)$ spaces originated as a generalization of Hadamard manifolds \cite{kiyoshi,bridson-haefliner-book}.
Most CAT$(0)$ spaces appearing in practice, such as positive definite matrices, and Hilbert spaces are complete.
A complete CAT$(0)$ space is called a Hadamard space. 
This is the class of spaces we are going to focus on. 
For a history, more applications, and classical references on CAT$(0)$ spaces we refer to the book of Bridson and Häfliger \cite{bridson-haefliner-book}. 
For a discussion of recent advances and open problems we refer to the survey article of Ba\v{c}\'{a}k \cite{bacak_challenges}, and for applications in optimization, we refer to the book of Ba\v{c}\'{a}k \cite{bacakbook}.

One of the fundamental problems in statistics is summarizing a given sample by its average.
Inspired by the Gaussian method of least squares \cite{gauss1809theoria}, Fr\'echet \cite{frechet-mean} proposed to consider the set of minimizers of
\begin{equation*}
  z \mapsto \frac{1}{n}\sum_{k=1}^{n}d(x_{k},z)^{2},
\end{equation*}
for a metric space $(N,d)$, and $z,x_{1},\dots,x_{n}\in N$.
While the definition of the Fr\'echet mean (also known as Karcher mean, or barycenter) is fairly straightforward, in a general metric space, the set of Fr\'echet means may be empty, or the mean of $n$ points may be non-unique \cite{ziezold1977, evans-jaffe-2024}.
In a Hadamard space the mean of $n$ points always exists and is unique \cite{bacakbook}.
However, even if Fr\'echet means exist and are unique, they are known to defy Euclidean intuition in surprisingly simple spaces such as open books \cite{hotz-sticky}.
As proper asymptotic and finite-sample guarantees for Fr\'echet means in Hadamard spaces have only been established recently \cite{le-gouic-fast-convergence,brunel,escande2023}, multiple alternative means have been proposed in the literature.
Es-Sahib and Heinich \cite{es-sahib-heinich} used an axiomatic approach to define a notion of mean on locally compact Hadamard spaces.
Navas \cite{navas-ergodic} was able to generalize the construction of Es-Sahib and Heinich to non-locally compact Hadamard spaces.
While these axiomatic means have nice theoretical properties, they are in general hard to compute in practice, as they frequently depend on some kind of limiting process.
On the other hand, computing the inductive mean introduced by Sturm \cite{sturm} is straightforward, as long as we can compute geodesics reasonably well. 
It is based on the observation that in Euclidean spaces, the mean of $n$ points is a convex combination of the first $n-1$ points and the $n$-th point.%, i.e.,
%\begin{equation*}
%  \frac{1}{n}\sum_{k=1}^{n}x_{k} = \biggl(1-\frac{1}{n}\biggr)\frac{1}{n-1}\sum_{k=1}^{n-1}x_{k} + \frac{1}{n}x_{n}.
%\end{equation*}
This idea can be generalized to Hadamard spaces and leads to a notion of mean, that is in general different from the Fr\'echet mean.
Given its importance, the case of symmetric, positive definite matrices has attracted considerable attention \cite{Hansen,holbrook-no-dice,bini2013,zhang2016,moakher2005}. 
Hansen \cite{Hansen} has introduced a notion of mean that is based on an idea similar to the inductive mean. 
Kim et al. \cite{kim} generalized this idea to the case of Hadamard spaces.
However, the theoretical properties of Hansen's mean are not as well understood as those of inductive means. 
In general, any two of these means may differ from one another (see Example 5.3 in \cite{kim}, and Example 6.5 in \cite{sturm}).

For most of these means, a more or less complete asymptotic theory is available.
Es-Sahib and Heinich \cite{es-sahib-heinich} provided a strong law of large numbers for their mean, and Navas \cite{navas-ergodic} was able to prove an $L^{1}$ ergodic theorem.
Sturm \cite{sturm} provides $L^{2}$ and strong laws of large numbers, together with an $L^{2}$ concentration inequality of order $O(n^{-1/2})$ for the inductive mean.
More recently, Antezana et al. \cite{antezana2023} proved an ergodic theorem for inductive means for $L^{1}$-functions from a compact, metrizable topological group, equipped with its Haar measure, into a Hadamard space.
Choi and Ji \cite{choi-toeplitz} provided a weighted version of the inductive mean and Hansen's mean and were able to give sufficient conditions for a strong law of large numbers.
For the Fr\'echet mean, early versions of a strong law of large numbers for general metric spaces were provided by Ziezold \cite{ziezold1977}, while the convergence of general $L^p$-\textit{centroids} was proved in \cite{slln-centroid-cpt-metric-space} for compact metric spaces.
Limit theorems for Fr\'echet means in general metric spaces have seen continued interest and progress \cite{bhattacharya-patrangenaru-2003, evans-jaffe-2024, schoetz}.
For CAT$(0)$ spaces, Le Gouic et al. \cite{le-gouic-fast-convergence} proved $O(n^{-1/2})$ concentration in $L^2$ of the Fr\'echet mean of an i.i.d. sample around the common population mean, given curvature bounds on the space. 
Recently, Brunel and Serres \cite{brunel} have shown that the distance between the population mean and both the inductive mean and the Fr\'{e}chet mean of points sampled independently from \textit{sub-Gaussian} distributions are of order $O(\log(1/\delta)n^{-1/2})$ with probability at least $1-\delta$, if we assume that the curvature of the underlying space is bounded from below.
This assumption can be dropped in the case of the Fr\'echet mean (see Section 6.1 in \cite{escande2023}). 
There neither sub-Gaussianity nor lower bounds on the curvature are required, rather the random variables are assumed to have finite $1$-exponential Orlicz norm.
Sch\"otz \cite{schötz18-entropy} and Ahidar-Coutrix et al. \cite{entropy-ahidarcoutrix} provided similar results under \textit{entropy conditions}.

All of these results have one thing in common: they require the underlying sample to be i.i.d.
In practice, any of these assumptions might be violated.
To confidently apply these results, practitioners require assurances regarding the stability of these laws. 
To the best of our knowledge, there are no results available in the literature, concerned with the stability of the theorem mentioned above.
This is precisely the aim of this paper.
Let us end this section by briefly outlining the contribution and structure of this note:
\begin{enumerate}
  \item In Theorem \ref{thm:slln-main} we provide strong and $L^{2}$ laws of large numbers for the inductive mean when the underlying data is independent, but not necessarily identically distributed.
  Additionally, we provide finite sample concentration inequalities in $L^2$.
    This generalizes a classical result of Sturm \cite{sturm}. 
    Brunel and Serres \cite{brunel} established a weak law of large numbers for independent, but not necessarily identically distributed, sub-Gaussian random variables with the same mean, as a consequence of a more general concentration inequality. 
    The $L^{p}$ and almost sure convergence of non-identically distributed random variables is still open.
    We are able to close this gap.
    Neither do we require our random variables to be sub-Gaussian, nor do we require them to have the same mean. 
    Additionally, the classical strong law of large numbers of Sturm \cite{sturm} requires the support of the random variables to be bounded, which we can relax considerably.
  \item In Theorem \ref{thm:limit-bn} we show that subject to reasonable regularity conditions on the first two moments of our sample, the only feasible $L^{2}$ limit of the Fr\'echet mean is the limit we recover in Theorem \ref{thm:slln-main}. 
    This can be regarded as a general $L^{2}$ version for Hadamard spaces of a classical theorem of Ziezold \cite{ziezold1977}, in which we neither require our data to be independent, nor identically distributed. 
  \item Proposition \ref{prop:stab-of-two-seq} contains a general stability result for both the Fr\'echet mean as well as the inductive mean. It neither requires the random variables to be independent nor identically distributed. 
 It is the basis for Propositions \ref{prob:huber-L2} and \ref{prop:huber-as}, in which we deal with a contamination model inspired by Huber's $\varepsilon$-contamination. In this model, we randomly replace our data sequence $X_{1},\dots,X_{n}$ with elements from a noise sequence $Y_{1},\dots,Y_{n}$. 
    In contrast to Huber's model, we neither require our data nor our noise to be identically distributed. In fact, Proposition \ref{prob:huber-L2} does not even require the noise to be independent of the data.
  \item In Section \ref{sec:means} we show that many means introduced in the literature satisfy a version of the Ces\`aro or Toeplitz Lemma. 
    Together with the previous results of Section \ref{sec:slln}, this means that roughly speaking, means in Hadamard spaces are no worse behaved than in linear spaces.
    Here we also introduce a stochastic resampling scheme, which we use to analyze the stability of the (non-symmetric) inductive mean with respect to information loss, sample permutation, and noise. 
  \item In Section 6 we perform two short simulation studies; one on the space of positive definite symmetric matrices, and one on the space of open books. 
    These studies highlight that the robustness of means depends dramatically on the geometry of the space. In the case of positive definite matrices a small degree of contamination is enough for various means to be non-consistent, while in the case of open books, some means converge, although $1/3$ of the data is contaminated.
\end{enumerate}

\section{Heteroscedastic Laws of Large Numbers}\label{sec:slln}

In this paper, we follow the analytic approach to Hadamard spaces popularized in the excellent article by Sturm \cite{sturm}.
For the classical viewpoint on Hadamard spaces via triangle comparison theorems, we recommend the book by Bridson and Häfliger \cite{bridson-haefliner-book}.
A CAT$(0)$ (or NPC space --  \textbf{N}on-\textbf{P}ositive-\textbf{C}urvature space) is a metric space $(H,d)$ such that for all $x,y\in H$, there is a $m\in H$, such that for all $z \in H$ 
  \begin{equation}\label{eq:char-cat(0)}
  d(z,m)^{2} \leq \frac{1}{2}d(z,x)^{2} + \frac{1}{2}d(z,y)^{2} - \frac{1}{4}d(x,y)^{2}.
\end{equation}
The point $m$ can be thought of as a midpoint between $x$ and $y$.
Intuitively speaking, triangles in CAT$(0)$ spaces are \textit{slimmer} than in Euclidean space.
  \begin{figure}[ht] 
    \begin{center}    
  \includegraphics[width=5cm]{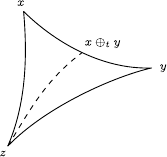}
    \end{center}
    \caption{Hyperbolicity of triangles in Hadamard spaces.}
  \end{figure}
A complete CAT$(0)$ space is called a Hadamard space. 
At first glance, this definition might seem somewhat opaque.
To fully appreciate this definition, we need some notions from metric geometry.
Let $(H,d)$ be a Hadamard space. 
A minimal geodesic between two points $x,y\in H$ is a continuous map $\gamma: [0,1] \to H$, with $\gamma(0)=x$, $\gamma(1) = y$, and 
\begin{equation*}
  d(\gamma(s),\gamma(t)) = |t-s| d(x,y) \quad \forall\, 0 \leq s \leq t \leq 1.
\end{equation*}
Minimal geodesics are of course geodesics, in the sense that
\begin{equation*}
  d(x,y) = \min_{\gamma\in C(x,y)}L(\gamma),
\end{equation*}
where $C(x,y)=\{\gamma:[0,1]\to H \mid \gamma\text{ continuous, } \gamma(0)=x, \gamma(1)=y\}$ is the set of all continuous curves from $[0,1]$ into $H$ connecting $x$ to $y$, and
\begin{equation*}
   L(\gamma) = \sup\biggl\{\sum_{i=0}^{n-1}d(\gamma(t_{i}),\gamma(t_{i+1})) \,\biggl\vert\, 0 = t_{0}< t_{1}< \cdots < t_{n} = 1\biggr\}
\end{equation*}
is the length of a curve $\gamma$ (the supremum is taken over all finite partitions of the interval $[0,1]$).
Inequality \eqref{eq:char-cat(0)} implies that any two points are connected by a unique minimal geodesic. 
In fact, Inequality \eqref{eq:char-cat(0)} \textit{extends} to minimal geodesics in the following sense.

\begin{proposition}[Proposition 2.3, \cite{sturm}]
  Let $(H,d)$ be a Hadamard space. Then any two points $x,y\in H$ are connected by a unique minimal geodesic $\gamma$. Furthermore, for all $z\in H$ and $t\in [0,1]$, we have
  \begin{equation}\label{eq:npc-ineq}
    d(z,\gamma(t))^{2} \leq (1-t)d(z,x)^{2} + td(z,y)^{2} - t(1-t)d(x,y)^{2}.
  \end{equation}
\end{proposition}
Inequality \eqref{eq:npc-ineq} is known as the NPC-inequality and is of fundamental importance in the theory of Hadamard spaces.
For any two points $x,y$ in a Hadamard space $H$, we write $x\oplus_{t}y \coloneqq \gamma(t)$ for the unique minimal geodesic $\gamma: [0,1] \to H$ joining $x$ and $y$. 
This allows us to introduce the inductive mean of Sturm \cite{sturm}. 
It is based on the following idea: in a normed space one may recursively compute the mean of a sequence $(x_{n})_{n\in \mathbb{N}}$ as  
\begin{equation*}
  \frac{1}{n}\sum_{k=1}^{n}x_{k} = \biggl(1-\frac{1}{n}\biggr)\frac{1}{n-1}\sum_{k=1}^{n-1}x_{k} + \frac{1}{n}x_{n}. 
\end{equation*}
In other words, the \textit{new} mean $n^{-1}\sum_{k=1}^{n}x_{k}$ lies on a geodesic between the \textit{old} mean $(n-1)^{-1}\sum_{k=1}^{n-1}x_{k}$ and the \textit{new} point $x_{n}$. 
This point of view can be generalized to Hadamard spaces. 
Let $H$ be a Hadamard space and $x_{1},\dots,x_{n}\in H$.
We set 
\begin{equation*}
  S_{1} = x_{1}, \quad \text{and} \quad S_{n+1} = S_{n}\oplus_{\frac{1}{n+1}} x_{n+1}.
\end{equation*}
It is worth pointing out, that $S_{n}$ can be computed by evaluating $n-1$ geodesics. 
In some sense, this is optimal, since we have to look at each of the $n$ points at least once to compute their Fr\'echet mean.
Furthermore, the inductive mean is contracting in quadratic mean, i.e., for all $z\in H$, we have
\begin{equation}\label{eq:sn-est}
  d(z,S_{n})^{2}\leq \frac{1}{n}\sum_{k=1}^{n}d(z,x_{k})^{2}.
\end{equation}
This is easily established by applying \eqref{eq:npc-ineq} iteratively and simply ignoring the negative terms.

However, the inductive mean may not be the most natural notion of mean on a Hadamard space, since it is based on the particular way one can compute the arithmetic mean in linear spaces (which is not unique, see Hansen's mean in Section \ref{sec:means}), rather than the classical least-squares optimization problem \cite{gauss1809theoria}.
This role is occupied by the Fr\'echet mean, which in general does not coincide with the inductive mean \cite{sturm}. 

To define the Fr\'echet mean, we may first talk about Hadamard space-valued random variables. 
A random variable on a Hadamard space $H$ is a Borel measurable function $X$ from some probability space into $H$. 
We say that $X$ is in $\mathcal{L}^{p}(H)$ for $p\geq 1$, if 
\begin{equation*}
  \mathbb{E}(d(X,z)^{p}) < \infty
\end{equation*}
for some (equivalently all) $z \in H$.
Clearly $\mathcal{L}^{p}(H) \subseteq \mathcal{L}^{q}(H)$ for $1 \leq p \leq q$.
We say $X_{n} \to Z$ in $\mathcal{L}^{p}(H)$, if $\mathbb{E}(d(X_{n},Z)^{p}) \to 0$.
For a real-valued random variable $X$ and $p\geq 1$, we denote with $\|X\|_{p} = \mathbb{E}(|X|^{p})^{1/p}$ its $L^{p}$-norm.
For some vector $x\in \mathbb{R}^n$, we denote with $\|x\|_{\ell^p} = \bigl(\sum_{k=1}^n |x_i|^p\bigr)^{1/p}$ its $\ell^p$-norm.
The expectation of a random variable $X\in \mathcal{L}^{2}(H)$ is given by
\begin{equation*}
  \mathbb{E}(X) = \argmin_{z\in H}\mathbb{E}(d(X,z)^{2}).
\end{equation*}
The $\argmin$ exists and is unique, by Proposition 2.2.17 in \cite{bacakbook}.
If $X\in \mathcal{L}^{1}(H)$, we can define its expectation by
\begin{equation*}
  \mathbb{E}(X)  = \argmin_{z\in H}\mathbb{E}(d(X,z)^{2} - d(X,y)^{2}),
\end{equation*}
for any $y\in H$.
Again, the $\argmin$ exists and is unique, since the function $z \mapsto \mathbb{E}(d(X,z)^{2} - d(X,y)^{2})$ is $1$-strongly convex (in the sense of \cite{bacakbook}) by the NPC inequality \eqref{eq:npc-ineq} and continuous. 
By Proposition 2.2.17 in \cite{bacakbook}, such functions have unique minimizers. 

We define the variance of $X\in \mathcal{L}^{2}(H)$ as $\mathrm{Var}(X) = \mathbb{E}(d(X,\mathbb{E}(X))^{2}) = \min_{z\in H}\mathbb{E}(d(X,z)^{2})$. 
The classical variance equality, $\mathrm{Var}(Y) = \mathbb{E}(Y^{2}) - \mathbb{E}(Y)^{2}$ for real-valued random variables $Y$, turns into a \textit{variance inequality} in Hadamard spaces.

\begin{lemma}[Variance Inequality, \cite{sturm}]\label{lemma:var-ineq}
For $X\in \mathcal{L}^{1}(H)$ and $z \in H$ we have
\begin{equation*}
  \mathbb{E}(d(X,z)^{2} - d(X,\mathbb{E}(X))^{2}) \geq d(z,\mathbb{E}(X))^{2}.
\end{equation*}
In particular, if $X\in \mathcal{L}^{2}(H)$, this can be written as
\begin{equation}\label{eq:var}
  \mathrm{Var}(X) \leq \mathbb{E}(d(X,z)^{2}) - \mathbb{E}(d(z,\mathbb{E}(X)))^{2}.
\end{equation}
\end{lemma}

Now we can define the Fr\'echet mean and discuss some of its properties. 
Let $x_{1},\dots,x_{n}$ be a sequence of points in a Hadamard space $H$.
The Fr\'echet mean (or barycenter) of this sequence is defined as
\begin{equation*}
  b_{n} = \argmin_{z\in H}\frac{1}{n}\sum_{k=1}^{n}d(x_{k},z)^{2}.
\end{equation*}
In other words, if $Y_{n}$ is a random variable with $\mathbb{P}(Y_n=x_{i}) = n^{-1}$, then $b_{n} = \mathbb{E}(Y_{n})$.
Applying Lemma \ref{lemma:var-ineq} to $Y_{n}$ yields
\begin{equation*}%\label{eq:var-ineq-barycenter}
  d(b_{n},z)^{2} \leq \frac{1}{n}\sum_{k=1}^{n}d(x_{k},z)^{2}-\frac{1}{n}\sum_{k=1}^{n}d(x_{k},b_{n})^{2}.
\end{equation*}

The Fr\'echet mean possesses several desirable theoretical properties. 
For example, the variance inequality (Lemma \ref{lemma:var-ineq}) is sharper than the $L^2$-contraction (Equation \eqref{eq:sn-est}) of the inductive mean.
Yokota generalized the classical Banach-Saks theorem to Fr\'echet means in Hadamard spaces \cite{yokota}.
Furthermore, it satisfies the canonical properties of a mean put forward by Ando et al. \cite{ando-li-mathias,Bhathia-Holbroo,lawson-lim-monotonic}.
However, computing the Fr\'echet mean is a non-trivial task.
We may look at the space of positive definite matrices as an example. 
Here explicit formulas for the distance and the unique minimal geodesic between two points are known. 
The distance $d_P$ between two positive definite matrices $A$ and $B$ is given by 
\begin{equation*}
  d_{P}(A,B) = \|\log(B^{-1/2}AB^{-1/2})\|_{F},
\end{equation*}
where $\|\cdot\|_{F}$ denotes the Frobenius norm.
The unique geodesic between $A$ and $B$ is given by 
\begin{equation*}
  t\mapsto A\oplus_{t}B = A^{1/2}(A^{-1/2}BA^{-1/2})^{t}A^{1/2} 
\end{equation*}
(see \cite{kim}).
This exemplifies a typical phenomenon in geodesic spaces -- in the worst case, one has to \textit{measure} the distance between two points by computing the length of a geodesic. 
In other words, computing a geodesic between two points is in general no more costly than computing the distance between them. 
Still, computing geodesics may be a non-trivial or computationally expensive task, as can be seen from the example above. 
Hence, one in general prefers to use an iterative approximation scheme, such as the proximal point algorithm (see Section \ref{sec:means}), to compute the Fr\'echet mean.
However, evaluating the objective function for the Fr\'echet mean once, requires the computation of $n$ distances, which in general is computationally as costly as computing the inductive mean.
%For the Fr\'echet mean, this implies that plugging in a single point in the optimization problem is about as expensive as computing the inductive mean.
In practice, data may also arrive in an online fashion, and one wants to update the predictors once new data is available. 
For the Fr\'echet mean, one has to solve a new optimization problem in general, while the inductive mean allows for easy online updates. 
Hence Sturm \cite{sturm} and subsequent authors \cite{antezana2023,choi-toeplitz} frequently focused on the inductive mean instead.   
However, it is an open problem whether $d(S_{n},b_{n}) \to 0$ for general sequences of points $(x_{n})_{n\in \mathbb{N}}$. 
%If the sequence is an i.i.d. sample drawn from a distribution with bounded support, a classical result of Sturm \cite{sturm} states that $S_{n}$ converges almost surely against the common mean of the underlying random sequence. 
%More recently, Brunel and Serres \cite{brunel} showed that the distance between the barycenter and the inductive mean $S_{n}$ of a sequence, independently sampled from \textit{sub-Gaussian} distributions, is of order $O(\log(\delta^{-1}) n^{-1/2})$ with probability at least $1-\delta$. 
%Escande \cite{escande2023} provided a similar result in $L^{2}$ for i.i.d. sub-exponential samples.
%Le Gouic et al. \cite{le-gouic-fast-convergence} showed that, subject to some geodesic extensibility conditions, the Fr\'echet mean of i.i.d. random variables with second moments converges in $L^{2}$ to the common mean of the sample.

While the $L^{2}$- and almost sure convergence of i.i.d. sequences is largely understood, the $L^{2}$- and almost sure convergence of non-identically distributed sequences and sequences sampled from distributions of unbounded support are still open. 
These questions are addressed in Theorem \ref{thm:slln-main} and Corollary \ref{cor:slln} respectively.
Theorem \ref{thm:slln-main} provides (strong) laws of large numbers for non-identically distributed but independent sequences of Hadamard space-valued random variables.
Corollary \ref{cor:slln} quantifies how fast the support of the random variables can grow, while a strong law of large numbers remains valid.

\begin{theorem}\label{thm:slln-main}
Let $X_n \in \mathcal{L}^2(H)$, $n\geq1$ be a sequence of independent random variables with expectation $\mathbb{E}(X_n) = \mu_n$, $\mu\in H$ and $D_{n}= \max_{1\leq k \leq n} \max\bigl\{d(\mu,\mu_{k}), \mathbb{E}(d(X_k,\mu_k))\bigr\}$. 
We have
\begin{equation}\label{eq:induction-claim}
  \mathbb{E}(d(S_n,\mu)^2) \leq \frac{6D_{n}}{n}\sum_{k=1}^{n}d(\mu_k,\mu) + \frac{1}{n^2} \sum_{k=1}^{n}\mathrm{Var}(X_k).
\end{equation}
%i.e., if 
%\begin{enumerate}
%  \item $D_{n}\frac{1}{n}\sum_{k=1}^{n} d(\mu_{k},\mu) \to 0,\quad$  and
%  \vspace{0.15cm}
%\item $\frac{1}{n^{2}} \sum_{k=1}^{n}\mathrm{Var}(X_{k}) \to 0$, 
%\end{enumerate}
%then $S_{n} \to \mu$ in $\mathcal{L}^{2}(H)$ and in probability. 
%If the sequence $(X_n)_{n\geq 1}$ is uniformly bounded almost surely (\textit{i.e.} all $X_{n}$ lie within $B_{r}(z)$ almost surely for, some $z\in H$ and $r>0$), and 
%  \begin{equation*}
%    \frac{1}{n}\sum_{k=1}^{n}d(\mu_{k},\mu) = O(n^{-p})
%  \end{equation*}
%  for some $p>\frac{1}{2}$, then, $S_{n} \to \mu$ almost surely.
\end{theorem}

%\begin{theorem}\label{thm:slln-main}
%  Let $X_n \in \mathcal{L}^2(H)$ be a sequence of independent random variables with expectation $\mathbb{E}(X_n) = \mu_n$, and $\mu\in H$ such that 
%\begin{enumerate}
%  \item $\sup_{n\in \mathbb{N}} d(\mu_{n},\mu) < \infty$,
%  \item $\frac{1}{n}\sum_{k=1}^{n} d(\mu_{k},\mu) \to 0$ and
%\item $\sup_{n\in \mathbb{N}} \mathrm{Var}(X_{n}) <\infty$. 
%\end{enumerate}
%Then $S_{n} \to \mu$ in $\mathcal{L}^{2}(H)$ and in probability. 
%  If the sequence $(X_n)_{n \in \mathbb{N}}$ is uniformly bounded almost surely and 
%  \begin{equation*}
%    \frac{1}{n}\sum_{k=1}^{n}d(\mu_{k},\mu) = O(k^{-p})
%  \end{equation*}
%  for some $p>\frac{1}{2}$, then, $S_{n} \to \mu$ almost surely.
%\end{theorem}
\begin{proof}
  %We are going to show by induction that
 % \begin{equation}
 % \mathbb{E}(d(S_n,\mu)^2) \leq \frac{6D_{n}}{n}\sum_{k=1}^{n}d(\mu_k,\mu) + \frac{1}{n^2} \sum_{k=1}^{n}\mathrm{Var}(X_k),
%\end{equation}
%which goes to zero by Assumptions 1 and 2.
We proceed by induction. For the case $n=1$, we have to show 
\begin{equation*}
  \mathbb{E}(d(X_{1},\mu)^{2}) \leq 6D_{1}d(\mu_{1},\mu) + \mathrm{Var}(X_{1}).
\end{equation*}
Since $\mathrm{Var}(X_{1}) = \mathbb{E}(d(X_{1},\mu_{1})^{2})$, we can rewrite this as
\begin{equation*}
\mathbb{E}(d(X_{1},\mu)^{2}-d(X_{1},\mu_{1})^2)  \leq 6D_{1} d(\mu_{1},\mu).
\end{equation*}
%Applying the Cauchy-Schwarz inequality, we get 
%\begin{equation*}
%  \begin{split} 
%    \mathbb{E}(d(X_{1},\mu)^{2}-d(X_{1},\mu_{1})^{2})  &= \mathbb{E}((d(X_{1},\mu)-d(X_{1},\mu_{1}))(d(X_{1},\mu)+d(X_{1},\mu_{1}))) \\
%                                                       & \leq \|d(X_{1},\mu) - d(X_{1},\mu_{1})\|_{2} \,\|d(X_{1},\mu) + d(X_{1},\mu_{1})\|_{2}.
%  \end{split}
%\end{equation*}
%By the reverse triangle inequality, we have 
%\begin{equation*}
%  \|d(X_{1},\mu) - d(X_{1},\mu_{1})\|_{2} \leq d(\mu_{1},\mu).
%\end{equation*}
By the reverse triangle inequality, we get 
\begin{equation}\label{eq:rev-triangle-argument}
  \begin{split} 
    \mathbb{E}(d(X_{1},\mu)^{2}-d(X_{1},\mu_{1})^{2})  &= \mathbb{E}((d(X_{1},\mu)-d(X_{1},\mu_{1}))(d(X_{1},\mu)+d(X_{1},\mu_{1}))) \\
                                                       & \leq d(\mu_{1},\mu)\|d(X_{1},\mu) + d(X_{1},\mu_{1})\|_{1}.
  \end{split}
\end{equation}
Furthermore, a simple application of the triangle inequality implies 
\begin{equation*}
  \|d(X_{1},\mu) + d(X_{1},\mu_{1})\|_{1} \leq 2\|d(X_{1},\mu_{1})\|_{1} + d(\mu_{1},\mu) \leq 3 D_{1}\leq 6D_{1},
\end{equation*}
proving the case $n=1$.

Moving on with the induction step, the NPC inequality \eqref{eq:npc-ineq} implies
\begin{equation}\label{eq:after-npc}
  \mathbb{E}(d(S_{n+1},\mu)^2) \leq \frac{n}{n+1}\mathbb{E}(d(S_n,\mu)^2) + \frac{1}{n+1}\mathbb{E}(d(\mu, X_{n+1})^2) - \frac{n}{(n+1)^2}\mathbb{E}(d(S_{n},X_{n+1})^2).
\end{equation}
We would like to apply the variance inequality (Lemma \ref{lemma:var-ineq}) to the last term. To this end we may write $\mathbb{E}(d(S_{n},X_{n+1})^2) = \mathbb{E}(\mathbb{E}(d(S_{n},X_{n+1})^2\mid S_n))$.
As $X_{n+1}$ and $S_n$ are independent, the expression $\mathbb{E}(d(X_{n+1},S_n)^2\mid S_n)$ may be written as $H(S_n)$ where $H(z) = \mathbb{E}(d(z,X_{n+1})^2)$, for $z\in H$. Now, by the variance inequality \eqref{eq:var}
\begin{equation*}
  H(z) \geq d(z,\mu_{n+1})^2 + \mathbb{E}(d(\mu_{n+1},X_{n+1})^2),
\end{equation*}
and hence 
\begin{equation*}
  \begin{split} 
    \mathbb{E}(d(S_{n},X_{n+1})^2) &= \mathbb{E}(\mathbb{E}(d(S_{n},X_{n+1})^2\mid S_n)) = \mathbb{E}(H(S_n)) \\
                                 & \geq \mathbb{E}(d(S_n,\mu_{n+1})^2) + \mathbb{E}(d(\mu_{n+1},X_{n+1})^2).
  \end{split}
\end{equation*}
Combining this with Inequality \eqref{eq:after-npc} yields
%\begin{equation*}
%  \begin{split} 
%    \mathbb{E}(d(S_{n+1},\mu)^2) &\leq \frac{n}{n+1}\mathbb{E}(d(S_n,\mu)^2) + \frac{1}{n+1}\mathbb{E}(d(\mu, X_{n+1})^2)\\
%                                 & -\frac{n}{(n+1)^2}\bigl\{\mathbb{E}(d(S_n,\mu_{n+1})^2) + \mathbb{E}(d(\mu_{n+1},X_{n+1})^2)\bigr\}.
%                                 %& =\frac{n}{n+1}\bigl[\mathbb{E}(d(S_{n},\mu)^{2}) - \frac{1}{n+1}\mathbb{E}(d(S_{n},\mu_{n+1})^{2})\bigr]\\
%                                 %&+ \frac{1}{n+1} \bigl[\mathbb{E}(d(\mu,X_{n+1})^{2}) - \frac{n}{n+1}\mathbb{E}(d(\mu_{n+1},X_{n+1})^{2})\bigr] = I+ II.
%  \end{split}
%\end{equation*}
%Regrouping terms gives
\begin{equation*}
  \mathbb{E}(d(S_{n+1},\mu)^2) \leq I + II,
\end{equation*}
where
\begin{equation*}
  \begin{split}
    I &= \frac{n}{n+1}\biggl\{\mathbb{E}(d(S_{n},\mu)^{2}) - \frac{1}{n+1}\mathbb{E}(d(S_{n},\mu_{n+1})^{2})\biggr\}, \quad \text{and}\\
    II & = \frac{1}{n+1} \biggl\{\mathbb{E}(d(\mu,X_{n+1})^{2}) - \frac{n}{n+1}\mathbb{E}(d(\mu_{n+1},X_{n+1})^{2})\biggr\}. 
  \end{split}
\end{equation*}
Starting with the term $I$, we insert $\pm n(n+1)^{-2}\mathbb{E}(d(S_{n},\mu)^{2})$. This yields 
\begin{equation*}
  I = \biggl(\frac{n}{n+1}\biggr)^{2} \mathbb{E}(d(S_{n},\mu)^{2}) + \frac{n}{(n+1)^{2}}\bigl\{\mathbb{E}(d(S_{n},\mu)^{2}) - \mathbb{E}(d(S_{n},\mu_{n+1})^{2})\bigr\}.
\end{equation*}
%Now looking at the second term of $I$, we may apply the Cauchy-Schwarz inequality to the effect of 
%\begin{equation*}
%  \begin{split}
%    \mathbb{E}(d(S_{n},\mu)^{2}) - \mathbb{E}(d(S_{n},\mu_{n+1})^{2}) & = \mathbb{E}((d(S_{n},\mu) - d(S_{n},\mu_{n+1}))(d(S_{n},\mu) + d(S_{n},\mu_{n+1}))) \\
%                                                                        &\leq \|d(S_{n},\mu) - d(S_{n},\mu_{n+1})\|_{2}\,\|d(S_{n},\mu) + d(S_{n},\mu_{n+1})\|_{2}.
%  \end{split}
%\end{equation*}
%By the reverse triangle inequality, we have 
%\begin{equation*}
%\|d(S_{n},\mu) - d(S_{n},\mu_{n+1})\|_{2} \leq d(\mu,\mu_{n+1}).
%\end{equation*}
For the second term of $I$, the same argument as in \eqref{eq:rev-triangle-argument} yields 
\begin{equation*}
  \begin{split}
    \mathbb{E}(d(S_{n},\mu)^{2}) - \mathbb{E}(d(S_{n},\mu_{n+1})^{2}) &\leq d(\mu,\mu_{n+1})\|d(S_{n},\mu) + d(S_{n},\mu_{n+1})\|_{1}.
  \end{split}
\end{equation*}
Applying the triangle inequality for the $L^{1}$-norm, we get 
\begin{equation}\label{eq:I-error-factor}
\|d(S_{n},\mu) + d(S_{n},\mu_{n+1})\|_{1} = \|d(S_{n},\mu)\|_{1} + \|d(S_{n},\mu_{n+1})\|_{1} \leq 2\|d(S_{n},\mu)\|_{1} + d(\mu,\mu_{n+1}).
\end{equation}
%Using the NPC inequality inductively on $d(S_{n},\mu)^{2}$, we arrive at the estimate
Using Proposition 3.1 in \cite{Funano}, $S_n$ is $L^1$-contracting, and we get 
\begin{equation*}
  d(S_{n},\mu) \leq \frac{1}{n}\sum_{k=1}^{n}d(X_{k},\mu).
\end{equation*}
%By the triangle inequality, and the estimate $(a+b)^{2}\leq 2(a^{2}+b^{2})$ for $a,b \geq 0$, we get
Now, the triangle inequality yields
\begin{equation*} 
  d(S_{n},\mu) \leq \frac{1}{n}\sum_{k=1}^{n}d(X_{k},\mu) \leq \frac{1}{n} \sum_{k=1}^{n}d(X_{k},\mu_{k}) +  \frac{1}{n}\sum_{k=1}^{n}d(\mu_{k},\mu).
\end{equation*}
Taking expectations, we get
\begin{equation*}
  \mathbb{E}(d(S_{n},\mu)) \leq \frac{1}{n}\sum_{k=1}^{n}\mathbb{E}(d(X_k,\mu_k)) + \frac{1}{n}\sum_{k=1}^{n}d(\mu_{k},\mu) \leq 2 D_n.
\end{equation*}
%and hence $ \| d(S_{n},\mu)\|_{1}\leq \|d(S_{n},\mu)\|_{2}\leq 2 D_{n}$ in \eqref{eq:I-error-factor}.
Since $D_{n}\leq D_{n+1}$, the right-hand side of Inequality \eqref{eq:I-error-factor} is bounded by $3D_{n+1}$, which implies 
\begin{equation*}
  \frac{n}{(n+1)^{2}}\bigl\{\mathbb{E}(d(S_{n},\mu)^{2}) - \mathbb{E}(d(S_{n},\mu_{n+1})^{2})\bigr\}\leq  \frac{3nD_{n+1}}{(n+1)^{2}} d(\mu,\mu_{n+1}) \leq \frac{3D_{n+1}}{n+1} d(\mu,\mu_{n+1}).
\end{equation*}
In total, this gives
\begin{equation*}%\label{eq:I-final}
  I \leq \biggl(\frac{n}{n+1}\biggr)^{2} \mathbb{E}(d(S_{n},\mu)^{2}) + \frac{3D_{n+1}}{n+1} d(\mu,\mu_{n+1}).
\end{equation*}
  For the second term $II$, we may insert $\pm (n+1)^{-1}\mathbb{E}(d(\mu_{n+1},X_{n+1})^{2})$. This yields
\begin{equation*}
  II = \frac{1}{(n+1)^{2}} \mathrm{Var}(X_{n+1}) + \frac{1}{n+1}\bigl\{\mathbb{E}(d(\mu,X_{n+1})^{2}) - \mathbb{E}(d(\mu_{n+1},X_{n+1})^{2})\bigr\}.
\end{equation*}
%By the Cauchy-Schwarz inequality, we may again estimate
%\begin{equation}\label{eq:II-CS}
%  \begin{split}
%    \mathbb{E}(d(\mu,X_{n+1})^{2}) - \mathbb{E}(d(\mu_{n+1},X_{n+1})^{2})& =\mathbb{E}((d(\mu,X_{n+1})-d(X_{n+1},\mu_{n+1}))(d(\mu,X_{n+1})+d(X_{n+1},\mu_{n+1})))\\
%                                                                         & \leq \|d(\mu,X_{n+1}) - d(X_{n+1},\mu_{n+1})\|_{2}\,\|d(X_{n+1},\mu) + d(X_{n+1},\mu_{n+1})\|_{2}.
%  \end{split}
%\end{equation}
%Applying the reverse triangle inequality again, we get
%\begin{equation*}
%  \|d(\mu,X_{n+1})-d(X_{n+1},\mu_{n+1})\|_{2} \leq d(\mu,\mu_{n+1}).
%\end{equation*}
Using the same idea as in \eqref{eq:rev-triangle-argument}, we get
\begin{equation}\label{eq:II-CS}
  \begin{split}
    \mathbb{E}(d(\mu,X_{n+1})^{2}) - \mathbb{E}(d(\mu_{n+1},X_{n+1})^{2}) \leq 3 d(\mu,\mu_{n+1})\,D_{n+1}.%\|d(X_{n+1},\mu) + d(X_{n+1},\mu_{n+1})\|_{1}.
  \end{split}
\end{equation}
%Applying the reverse triangle inequality again, we get
%\begin{equation*}
%  \|d(\mu,X_{n+1})-d(X_{n+1},\mu_{n+1})\|_{2} \leq d(\mu,\mu_{n+1}).
%\end{equation*}
%To estimate the second factor of the right-hand side of Inequality \eqref{eq:II-CS}, we may observe that $\mathbb{E}(d(X_{n+1},\mu_{n+1})^{2}) \leq \mathbb{E}(d(X_{n+1},\mu)^{2})$, as $\mathrm{Var}(X_{n+1}) = \min_{z\in H}\mathbb{E}(d(X_{n+1},z)^{2}) = \mathbb{E}(d(X_{n+1},\mu_{n+1})^{2})$. Hence a rough estimate yields
%A simple application of the triangle inequality yields
%\begin{equation*}
%  \begin{split} 
%    \|d(\mu,X_{n+1})+d(X_{n+1},\mu_{n+1})\|_{1} &\leq 2\mathbb{E}(d(X_{n+1},\mu_{n+1})) + d(\mu,\mu_{n+1}) \leq 3D_{n+1}.
%  \end{split}
%\end{equation*}
In total, this gives
\begin{equation*}
  II \leq \frac{1}{(n+1)^{2}} \mathrm{Var}(X_{n+1}) + \frac{3D_{n+1}}{n+1}d(\mu,\mu_{n+1}).
\end{equation*}
Combining the estimates for $I$ and $II$, and using the induction hypothesis as well as the monotonicity of $D_{n}$ in $n$, yields
\begin{equation*}
  \begin{split}
    \mathbb{E}(d(S_{n+1},\mu)^{2}) &\leq \biggl(\frac{n}{n+1}\biggr)^{2}\mathbb{E}(d(S_{n},\mu)^{2}) + \frac{3D_{n+1}}{n+1}d(\mu,\mu_{n+1}) + \frac{1}{(n+1)^{2}} \mathrm{Var}(X_{n+1}) + \frac{3D_{n+1}}{n+1}d(\mu,\mu_{n+1}) \\
                                   & \leq \frac{6D_{n+1}}{n+1}\sum_{k=1}^{n+1}d(\mu,\mu_{k}) + \frac{1}{(n+1)^{2}}\sum_{k=1}^{n+1}\mathrm{Var}(X_{k}).
  \end{split}
\end{equation*}
%Using the induction hypothesis and the monotonicity of $D_{n}$ in $n$, we get 
%\begin{equation*}
%  \mathbb{E}(d(S_{n+1},\mu)^{2}) \leq \frac{6D_{n+1}}{n+1}\sum_{k=1}^{n+1}d(\mu,\mu_{k}) + \frac{1}{(n+1)^{2}}\sum_{k=1}^{n+1}\mathrm{Var}(X_{k}).
%\end{equation*}
%This proves \eqref{eq:induction-claim}, and implies that $S_{n}\to \mu$ in $\mathcal{L}^{2}(H)$ and in probability.

\end{proof}

\begin{corollary}
Let $X_n \in \mathcal{L}^2(H)$, $n\geq1$ be a sequence of independent random variables with expectation $\mathbb{E}(X_n) = \mu_n$, $\mu\in H$ and $D_{n}= \max_{1\leq k \leq n} \max\bigl\{d(\mu,\mu_{k}), \mathbb{E}(d(X_k,\mu_k))\bigr\}$. 
If 
\begin{enumerate}
  \item $D_{n}\frac{1}{n}\sum_{k=1}^{n} d(\mu_{k},\mu) \to 0,\quad$  and
  \vspace{0.15cm}
\item $\frac{1}{n^{2}} \sum_{k=1}^{n}\mathrm{Var}(X_{k}) \to 0$, 
\end{enumerate}
then $S_{n} \to \mu$ in $\mathcal{L}^{2}(H)$ and in probability. 
\end{corollary}
If $\mu_{n} = \mu$ for all $n$ we recover Lemma 4.2 of Brunel and Serres \cite{brunel}.
In particular, it is worth pointing out, that we do not assume that our random variables have the same mean, nor that their means converge to some $\mu$. 
The term $D_n$ appearing in Theorem \ref{thm:slln-main} can be thought of as a maximal standard deviation. 
In this sense, the first condition of Theorem \ref{thm:slln-main} essentially means that $\mu_{n}$ converges to $\mu$ in Ceasaro mean faster than the standard deviation diverges. 
Intuitively speaking, the \textit{signal beats the noise}. To the best of our knowledge, these are the weakest assumptions in the current literature imposed on $\mu_{n}$ and $\mu$ which still guarantee laws of large numbers.
The second assumption requires the average variance to grow sublinearly.
This is similar, although slightly stronger than classical laws of large numbers require (cf. Theorem 2.3.10 in \cite{sen-singer}). 
In particular, if our random variables are i.i.d., we have a $O(n^{-1/2})$ concentration around $\mu$ in $L^2(H)$, similar to \cite{sturm, le-gouic-fast-convergence}.

In contrast to Corollary 1, the next Corollary is concerned with the almost sure convergence of $S_{n}$.
We may allow the support of the $X_{n}$'s to grow as $n$ goes to infinity. 
In this sense, Corollary \ref{cor:slln}, can be thought of as a generalization of Theorem 4.7 in \cite{sturm}.

%A closer inspection of our argument shows, that we do not have to require our $X_{n}$'s to be uniformly bounded for a strong law of large numbers to hold. 
%To be more precise, if $X_{1},\dots, X_{n}$ all lie within a ball whose radius grows \textit{slowly} with $n$, our strong law of large numbers still applies. 
%This is summarized in the following corollary.

\begin{corollary}\label{cor:slln}
Let $(X_{n})_{n\geq1}$ be a sequence of independent, Hadamard space-valued random variables with $\mathbb{E}(X_{n}) = \mu_{n}$, and let $\mu\in H$, such that 
\begin{equation*}
  \frac{1}{n}\sum_{k=1}^{n}d(\mu_{k},\mu) = O(n^{-p})
\end{equation*}
for some $p>\frac{1}{2}$. If there exist $z\in H$, $C>0$, and $0\leq q <\min \bigl\{\frac{1}{4},p-\frac{1}{2}\bigr\}$, such that 
\begin{equation}\label{eq:cor:supp-bound}
  \mathbb{P}\biggl(\max_{1\leq k\leq n}d(X_{k},z) \leq Cn^{q}\biggr) = 1,
\end{equation}
then $S_{n} \to \mu$ almost surely.
\end{corollary}

\begin{proof}
We follow the approach of Sturm \cite{sturm}. 
Note that \eqref{eq:cor:supp-bound} implies
\begin{equation*}
   \mathbb{P}\biggl(\max_{1\leq k\leq n}d(X_{k},\mu) \leq C'n^{q}\biggr) = 1,
\end{equation*}
for some $C'>0$.
We are going to show that $D_{n} = 2C'n^{q}$.
Since $\mu_{k}$ is the minimizer of $x\mapsto \mathbb{E}(d(X_{k},x)^{2})$, we have 
\begin{equation*}
  \mathbb{E}(d(X_{k},\mu_{k}))\leq \sqrt{\mathbb{E}(d(X_{k},\mu_{k})^{2})} \leq \sqrt{\mathbb{E}(d(X_{k},\mu)^{2})} \leq  C'k^{q},
\end{equation*}
which immedially implies
\begin{equation*}
  d(\mu,\mu_{k}) \leq \mathbb{E}(d(\mu,X_{k}))+ \mathbb{E}(d(X_{k},\mu_{k})) \leq 2C'k^{q}.
\end{equation*}
Next, we are going to show that $S_{n^{2}} \to \mu$.
%Since the $X_{n}$'s are almost surely uniformly bounded, their means $\mu_{n}$ and variances $\mathrm{Var}(X_{n})$, and thus $D_{n}$ are bounded as well.
%In particular, Assumptions 1 and 2 of the first part of the theorem are met. 
Inequality \eqref{eq:induction-claim} yields for every $\varepsilon>0$
\begin{equation*}
  \begin{split}
    \sum_{k=1}^{\infty} \mathbb{P}(d(S_{k^{2}},\mu)>\varepsilon) & \leq \sum_{k=1}^{\infty} \frac{1}{\varepsilon^{2}}\mathbb{E}(d(S_{k^{2}},\mu)^{2}) \\
                                                                 & \leq \frac{6}{\varepsilon^{2}} \sum_{k=1}^{\infty} \frac{D_{k^{2}}}{k^{2}}\sum_{j=1}^{k^{2}}d(\mu_{j},\mu) + \frac{1}{\varepsilon^{2}}\sum_{k=1}^{\infty} \frac{1}{k^{4}} \sum_{j=1}^{k^{2}}\mathrm{Var}(X_{j}) = I + II.
  \end{split}
\end{equation*}
Since $D_{k^{2}} \leq 2C'k^{2q}$, and $q< p-1/2$, the first term can be estimated by  
\begin{equation*}
  I \leq C_{0} \sum_{k=}^{\infty} k^{2(q-p)} <\infty,
\end{equation*}
for some constant $C_{0}>0$.
For the second term we observe that $\mathrm{Var}(X_{j})\leq C''j^{2q}$ for some $C''>0$.
Since $q<1/4$, this implies 
\begin{equation*}
  II \leq \frac{C''}{\varepsilon^{2}}\sum_{k=1}^{\infty}\frac{1}{k^{4}}\sum_{j=1}^{k^{2}}j^{2q} \leq \frac{C''_{0}}{\varepsilon^{2}}\sum_{k=1}^{\infty} k^{4(q-1)+2}< \infty,
\end{equation*}
for some constant $C_{0}''>0$.
%where  $C>0$, such that $6\sup_{n\in \mathbb{N}}D_{n} \leq C$.
Then, by the Borel-Cantelli lemma, $S_{n^{2}} \to \mu$ almost surely.  
%Since the sequence $(X_n)_{n\in\mathbb{N}}$ is uniformly bounded almost surely, there is some $z\in H$ and $r>0$ such that $d(X_{n},z) \leq r$ for all $n$ almost surely. 
Assumption \eqref{eq:cor:supp-bound} implies that $d(X_{n},\mu) \leq Cn^{q}$ almost surely.
Inequality \eqref{eq:sn-est} implies that $d(S_{n},\mu) \leq Cn^{q}$ for all $n$ almost surely. Hence we have 
\begin{equation*}
  d(S_{n},S_{n+1}) \leq \frac{1}{n+1}d(S_{n},X_{n+1}) \leq \frac{2n^{q}}{n+1},
\end{equation*}
for all $n$ almost surely. Thus we almost surely have for all $n$ and $n^{2} \leq k < (n+1)^{2}$
\begin{equation*}
  d(S_{n^{2}},S_{k}) \leq 2C \sum_{l=n^{2}}^{k} \frac{l^{q}}{l+1} \leq 2C \frac{(2n+1)(n+1)^{2q}}{n^{2}} = O(n^{-1/2}),
\end{equation*}
since $q\leq 1/4$.
This implies the strong law of large numbers. 

\end{proof}
%If the sequence $(X_n)_{n\geq 1}$ is uniformly bounded almost surely (\textit{i.e.} all $X_{n}$ lie within $B_{r}(z)$ almost surely for, some $z\in H$ and $r>0$), and 
%  \begin{equation*}
%    \frac{1}{n}\sum_{k=1}^{n}d(\mu_{k},\mu) = O(n^{-p})
%  \end{equation*}
%  for some $p>\frac{1}{2}$, then, $S_{n} \to \mu$ almost surely.

In a classical paper, Ziezold \cite{ziezold1977} showed that in a separable quasi-metric space, almost any limit point of the sequence of Fr\'echet means $b_{n}$ of an i.i.d. sequence $X_{1},\dots,X_{n}$ is a minimizer of $z\mapsto \mathbb{E}(d(X_{1},z)^{2})$.
However, if the sequence $X_{1},\dots,X_{n}$ is not i.i.d., one can not expect limit points to be minimizers of $z\mapsto \mathbb{E}(d(X_1,z)^2)$.
For Hadamard spaces, Theorem \ref{thm:limit-bn} below characterizes any potential $L^{2}$ limit point of the sequence of Fr\'echet means in a non-i.i.d. setting. 
We show that any such limit point must be the same limit we recover in Theorem \ref{thm:slln-main} using the inductive mean, i.e., $\mu$.
As such, Theorem \ref{thm:limit-bn} can be thought of as an $L^{2}$-version of Ziezold's result on Hadamard spaces.
It is worth pointing out, that we neither require our data to be independent, nor identically distributed. 
Additionally, we do not impose curvature bounds on the underlying space, nor do we impose geodesic extensibility conditions, as is the case in \cite{le-gouic-fast-convergence, brunel}.
Imposing some reasonably mild conditions on the first two moments of our sample is enough. 
However, whether the conditions below ensure convergence in general is open.

  \begin{theorem}\label{thm:limit-bn}
  Assume $\sup_{n\in \mathbb{N}}\mathrm{Var}(X_{n})<\infty$, and $n^{-1}\sum_{k=1}^{n}d(\mu,\mu_{k})^{2}\to 0$. If a subsequence $b_{n_{m}}$ of $b_{n}$ converges to some point $z\in H $ in $\mathcal{L}^{2}(H)$, then $z=\mu$. 
%  In particular, if $X_{1},\dots,X_{n}$ belong to a convex and %compact subset of $\mathcal{L}^{2}(H)$, $b_{n}$ converges to $\mu$ in $\mathcal{L}^{2}(H)%(H)$.
\end{theorem}

%{\color{red} Achtung Fehler: Man weiß im zweiten Teil nur, dass %$b_{n}$ gegen ein $Y\in \mathcal{L}^{2}(H)$ konvergiert. Damit man den %ersten Teil anwenden kann, muss aber $Y$ a.s. konstant sein. }
\begin{proof}
Let $b_{n_{m}}$ be a subsequence of $b_{n}$ converging to some $z\in H$ in $\mathcal{L}^{2}(H)$. Let us define the function
\begin{equation*}
  G(a) = \limsup_{m\to \infty}\frac{1}{n_{m}}\sum_{k=1}^{n_{m}}\mathbb{E}(d(X_{k},a)^{2}),
\end{equation*}
for $a \in H$. First, we are going to show, that $G$ has a unique minimizer.
Then, we are going to show, that any $\mathcal{L}^{2}(H)$ limit $z$ of $b_{n_m}$ is a minimizer of $G$, and finally, we are going to show that $G(\mu)\leq G(z)$, which implies $z=\mu$.

%Since $\sqrt{x}$ is monotonically increasing, and continuous, we have 
%\begin{equation*}
%  G(a) = \sqrt{\liminf_{n\to \infty}\frac{1}{n} \sum_{k=1}^{n}\mathbb{E}(d(X_{k},a)^{2})}.
%\end{equation*}
Let $a_{t}$ be the unique, minimal geodesic between $a_{0}$ and $a_{1}$. The function $G$ satisfies
\begin{equation*}
G(a_{t}) \leq (1-t)G(a_{0}) + tG(a_{1}) - t(1-t)d(a_{0},a_{1})^2.
\end{equation*}
In other words, $G$ is a strongly convex function in the sense of \cite{bacakbook}.
By Proposition 2.2.17 in \cite{bacakbook}, $G$ has a unique minimizer.

Now, we are going to show that every $\mathcal{L}^{2}(H)$ limit $z$ of $b_{n_{m}}$ is a minimizer of $G$. Using the triangle inequality, we have
\begin{equation}\label{eq:bn->z=>z=min}
  \begin{split}
    \frac{1}{n_{m}}\sum_{k=1}^{n_{m}}\mathbb{E}(d(X_{k},z)^{2}) & \leq \frac{1}{n_{m}}\sum_{k=1}^{n_{m}}\mathbb{E}(d(X_{k},b_{n_{m}})^{2}) + \frac{2}{n_{m}}\sum_{k=1}^{n_{m}}\mathbb{E}(d(b_{n_{m}},z)d(X_{k},b_{n_{m}})) + \frac{1}{n_{m}}\sum_{k=1}^{n_{m}}\mathbb{E}(d(b_{n_{m}},z)^{2}).
  \end{split}
\end{equation}
Since $\mathbb{E}(d(b_{n_{m}},z)^2)\to 0$ by assumption, the last term converges to zero. 
For the second term, Hölder's inequality implies
\begin{equation*}
  \frac{2}{n_{m}}\sum_{k=1}^{n_{m}}\mathbb{E}(d(b_{n_{m}},z)d(X_{k},b_{n_{m}}))\leq \frac{2}{n_{m}}\sum_{k=1}^{n_{m}}\sqrt{\mathbb{E}(d(b_{n_{m}},z)^{2})}\sqrt{\mathbb{E}(d(X_{k},b_{n_{m}})^{2})}.
\end{equation*}
Since 
\begin{equation*}
  \mathbb{E}(d(b_{n_{m}},X_{k})^{2})\leq \frac{1}{n_{m}}\sum_{l=1}^{n_{m}}\mathbb{E}(d(X_{l},X_{k})^{2}) \leq \frac{4}{n_{m}}\sum_{l=1}^{n_{m}}\mathbb{E}(d(X_{l},\mu_{l})^{2}) + d(\mu_{l},\mu)^{2} + d(\mu,\mu_{k})^{2} + \mathbb{E}(d(\mu_{k},X_{k})^{2}),
\end{equation*}
where the first inequality is a consequence of the variance inequality \ref{lemma:var-ineq}, and the second follows by applying the triangle inequality, and $\|x\|_{\ell^{1}}^{2}\leq 4 \|x\|_{\ell^{2}}^{2}$ for $x \in \mathbb{R}^{4}$. The last expression is uniformly bounded in $m$ by our assumptions.
This, together with $\mathbb{E}(d(b_{n_{m}},z)^{2})\to 0$, implies that the second term in \eqref{eq:bn->z=>z=min} vanishes.
In other words,
\begin{equation*}
  \limsup_{m\to \infty}\frac{1}{n_{m}}\sum_{k=1}^{n_{m}}\mathbb{E}(d(X_{k},z)^{2})\leq \limsup_{m\to \infty}\frac{1}{n_{m}}\sum_{k=1}^{n_{m}}\mathbb{E}(d(X_{k},b_{n_{m}})^{2}).
\end{equation*}
However, for each $\omega\in \Omega$, we have
\begin{equation*}
  b_{n_{m}}(\omega) = \argmin_{a\in H} \frac{1}{n_{m}}\sum_{k=1}^{n_{m}}d(X_{k}(\omega),a)^{2},
\end{equation*}
which in particular implies
\begin{equation*}
 G(z) \leq \limsup_{m\to \infty}\frac{1}{n_{m}}\sum_{k=1}^{n_{m}}\mathbb{E}(d(X_{k},b_{n_{m}})^{2}) \leq \limsup_{m\to \infty}\frac{1}{n_{m}}\sum_{k=1}^{n_{m}}\mathbb{E}(d(X_{k},a)^{2})  = G(a),
\end{equation*}
for every $a\in H$, i.e, $z$ is a minimizer of $G$.

Next, we are going to show that $G(\mu)\leq G(z)$. 
The triangle inequality implies
\begin{equation}\label{eq:G(m)<G(z)-est}
  \begin{split}
    G(\mu) & = \limsup_{m\to \infty}\frac{1}{n_{m}}\sum_{k=1}^{n_{m}}\mathbb{E}(d(X_{k},\mu)^{2}) \leq \limsup_{m\to \infty}\frac{1}{n_{m}}\sum_{k=1}^{n_{m}}\mathbb{E}(d(X_{k},\mu_{k})^{2})\\
           & + \limsup_{m\to \infty}\frac{2}{n_{m}}\sum_{k=1}^{n_{m}}d(\mu,\mu_{k})\mathbb{E}(d(X_{k},\mu_{k})) + \limsup_{m\to \infty}\frac{1}{n_{m}}\sum_{k=1}^{n_{m}}d(\mu,\mu_{k})^{2}.
  \end{split}
\end{equation}
The last term converges to zero, by our assumptions.
For the second term, Hölder's inequality implies
\begin{equation*}
  \begin{split}
    \frac{1}{n_{m}}\sum_{k=1}^{n_{m}}d(\mu,\mu_{k})\mathbb{E}(d(X_{k},\mu_{k})) &\leq \sqrt{\frac{1}{n_{m}}\sum_{k=1}^{n_{m}}d(\mu,\mu_{k})^{2}}\sqrt{\frac{1}{n_{m}}\sum_{k=1}^{n_{m}}\mathbb{E}(d(X_{k},\mu_{k}))^{2}} \\
                                                                                &\leq \sup_{k\in \mathbb{N}}\mathrm{Var}(X_{k}) \sqrt{\frac{1}{n_{m}}\sum_{k=1}^{n_{m}}d(\mu,\mu_{k})^{2}} \to 0. 
  \end{split}
\end{equation*}
For the first term, we recall that 
\begin{equation*}
  \mu_{k}= \argmin_{a\in H} \mathbb{E}(d(X_{k},a)^{2}),
\end{equation*}
and hence
\begin{equation*}
  \limsup_{m\to \infty}\frac{1}{n_{m}}\sum_{k=1}^{n_{m}}\mathbb{E}(d(X_{k},\mu_{k})^{2}) \leq \limsup_{m\to \infty}\frac{1}{n_{m}}\sum_{k=1}^{n_{m}}\mathbb{E}(d(X_{k},z)^{2}) = G(z).
\end{equation*} 
Summing up, \eqref{eq:G(m)<G(z)-est} implies that $\mu$ is a minimizer of $G$.

\end{proof}

%\begin{proposition}
%Let $X_{1},\dots,X_{n}\in \mathcal{L}^{2}(H)$ be i.i.d. with mean $\mu$, let $r_{n}\leq n$ be the contamination ratio, and let $Y_{1},\dots,Y_{r_{n}}\in \mathcal{L}^{2}(H)$ i.i.d. be the contamination.
%We assume that $\{X_{1},\dots,X_{n}, Y_{1},\dots,Y_{r_{n}}\}$ is independent. 
%The sequence $Z_{1},\dots,Z_{n}$ is generated by replacing $r_{n}$ of the $X_{i}$'s with $Y_{j}$'s (avoiding duplication). 
%Let $S_{n}$ be the inductive mean of $Z_{1},\dots,Z_{n}$.
%If $r_{n} = o(n)$, $S_{n}\to \mu$ in $\mathcal{L}^{2}(H)$. 
%\end{proposition}
%
%
%\begin{proof}
%  We first deal with the $\mathcal{L}^{2}(H)$ convergence of $S_{n}$. 
%  To this end, it is enough to checkt that $Z_{1},\dots,Z_{n}$ satisfy the assumptions of the first part of Theorem \ref{thm:slln-main}. 
%Under the present assumptions, $D_{n}$ in Theorem \ref{thm:slln-main} is equal to $D_{n}= \max\{\mathrm{Var}(X_{1}), \mathrm{Var}(Y_{1})\} =D$, which is independent of $n$. 
%Hence, if we denote the mean of $Y_{1}$ with $\nu$, and the mean of $Z_{k}$ with $\mu_{k}$, we get
%\begin{equation*}
%  \frac{D_{n}}{n}\sum_{k=1}^{n}d(\mu_{k},\mu) = D d(\mu,\nu)\frac{r_{n}}{n} \to 0.
%\end{equation*}
%Since the variances of $X_{1}$ and $Y_{1}$ are finite, the second assumption of the first half of Theorem \ref{thm:slln-main} is automatically fulfilled. 
%\end{proof}

Proposition \ref{prop:stab-of-two-seq}--\ref{prob:huber-L2} deal with contaminated data.
Our contamination model is inspired by Huber's $\varepsilon$-contamination model \cite{huber64,huber65}, in which data is drawn from a parametric model $G_{\theta,\varepsilon} = (1-\varepsilon)F_\theta + \varepsilon C$, where $F_\theta$ is the CDF of our data, and $C$ is the CDF of some noise distribution.
This is somewhat impractical for general Hadamard spaces. 
In particular, we neither want our data, nor our noise to be identically distributed. 
Hence, we use a random replacement model, in which a certain amount $1 \leq r_n\leq n$ of our data $X_1,\dots,X_n$ is randomly replaced with noise $Y_1,\dots,Y_{r_n}$.
The ratio $r_n/n$ roughly corresponds to the $\varepsilon$ in Huber's $\varepsilon$-contamination model.
If $Z_1,\dots,Z_n$ is the sequence that is generated in that matter, we are concerned with the following question: Under which conditions do the inductive and Fr\'echet mean of the $Z_i$'s and $X_i$'s converge to the same limit?
Proposition \ref{prob:huber-L2} is based on Proposition \ref{prop:stab-of-two-seq}, and deals with the $L^2$-stability of the law of large numbers under contamination.
Proposition \ref{prop:huber-as} is based on Corollary \ref{cor:slln}, and deals with the almost sure case. 

It is worth pointing out that the replacement scheme for the $L^2$ case (Proposition \ref{prob:huber-L2}) and the almost sure case (Proposition \ref{prop:huber-as}) differ slightly. 
In Proposition \ref{prob:huber-L2}, we fix an $r_n$ between $1$ and $n$, and pick a random function $f$, uniformly from
  \begin{equation*}
    \mathcal{F}=\{f:\{1,\dots,n\} \to \{0,1\}\mid |f^{-1}(1)| = r_{n}\},
\end{equation*}
where $|A|$ denotes the cardinality of a set $A$.
Note that $\mathcal{F}$ has finite cardinality. 
We use this random function $f$ to decide which data points are replaced with noise, i.e.,
\begin{equation*}
    Z_{k} = \begin{cases}
      X_{k} & \text{if } f(k) = 0 \\
      Y_{k} & \text{if } f(k) = 1
    \end{cases},
\end{equation*}
if $X_1,\dots,X_n$ is our data sequence, and $Y_1,\dots,Y_n$ is our noise sequence. 
On the other hand, in Proposition \ref{prop:huber-as} we replace a fixed (deterministic) but unknown number $r_n$ of our sample with noise. 
Random replacement is possible, in this case. 
However, we believe the setup used in Proposition \ref{prop:huber-as} to be more natural.

%where the contamination comes from a noise sequence and contaminates the data with a certain rate, which depends on $n$. Proposition Proposition \ref{prop:stab-of-two-seq} serves as basis for Proposition \ref{prop:huber-as} and \ref{prob:huber-L2}

  \begin{proposition}\label{prop:stab-of-two-seq}
Let $(X_{k})_{k\geq 1}$ and $(Y_{k})_{k\geq 1}$ be two sequences of Hadamard space-valued random variables.
Denote with $S_{n}^{X}$ and $b_{n}^{X}$ the inductive and Fr\'echet mean of $X_{1},\dots, X_{n}$, and with $S_{n}^{Y}$ and $b_{n}^{Y}$ the inductive and Fr\'echet mean of $Y_{1},\dots, Y_{n}$ respectively.
%let $S_{n}^{X}$ and $S_{n}^{Y}$ be the inductive means of $X_{1},\dots,X_{n}$ and $Y_{1},\dots, Y_{n}$ respectively. 
If $S_{n}^{X}$ converges in $\mathcal{L}^{2}(H)$ to $\mu$, and 
\begin{equation}\label{eq:SX-vs-SY}
  \frac{1}{n}\sum_{k=1}^{n}\mathbb{E}(d(X_{k},Y_{k})^{2}) \to 0,
\end{equation}
then $S_{n}^{Y} \to \mu$ in $\mathcal{L}^{2}(H)$. If in addition to \eqref{eq:SX-vs-SY}, $\mathbb{E}(\max_{1\leq k\leq n}d(X_{k},a)^{2}) = o(n)$ for some (and hence all) $a\in H$, and $b_{n}^{X}\to \mu$ in $\mathcal{L}^{2}(H)$, then $b_{n}^{Y} \to \mu$ in $\mathcal{L}^{2}(H)$. 
\end{proposition}

\begin{proof}
  If $x_{t}$ is the minimal geodesic between $x_{0}$ and $x_{1}$ and $y_{t}$ is the minimal geodesic between $y_{0}$ and $y_{1}$, then Lemma 1.2.2 in \cite{bacakbook} implies
  \begin{equation*}
    \begin{split}
      d(x_{t},y_{t})^{2} &\leq (1-t)^{2}d(x_{0},y_{0})^{2}+ t^{2}d(x_{1},y_{1})^{2} + t(1-t) \bigl[d(x_{0},y_{1})^{2} + d(x_{1},y_{0})^{2} - d(x_{0},x_{1})^{2}- d(y_{0},y_{1})^{2}\bigr].
    \end{split}
  \end{equation*}
  An application of Lemma 1.2.5 in \cite{bacakbook} to the factor of $t(1-t)$ in the last display yields
  \begin{equation}\label{eq:d(xt,yt)}
    \begin{split}
      d(x_{t},y_{t})^{2}& \leq (1-t)^{2}d(x_{0},y_{0})^{2} + t^{2}d(x_{1}y_{1})^{2} + t(1-t) \bigl[d(x_{0},y_{0})^{2} + d(x_{1},y_{1})^{2}\bigr] \\ 
                        & = (1-t) d(x_{0},y_{0})^{2} + td(x_{1},y_{1})^{2}.
    \end{split}
  \end{equation}
  Applying the last estimate inductively, we get
  \begin{equation*}
    \mathbb{E}(d(S_{n}^{X},S_{n}^{Y})^{2})\leq \frac{n-1}{n}\mathbb{E}(d(S_{n-1}^{X},S_{n-1}^{Y})^{2}) + \frac{1}{n}\mathbb{E}(d(X_{n},Y_{n})^{2}) \leq \cdots \leq \frac{1}{n}\sum_{k=1}^{n}\mathbb{E}(d(X_{k},Y_{k})^{2}),
  \end{equation*}
  and hence $S_{n}^{Y}\to \mu$ in $\mathcal{L}^{2}(H)$ by \eqref{eq:SX-vs-SY}.

  Moving on to the case of Fr\'echet means, we quickly recall a special case of the approximation scheme for $b_{n}^{X}$ introduced by Lim and P\'alfia \cite{lim-palfia}.
  There $b_{n}^{X}$ is approximated by an inductive mean based on an auxiliary sequence $Z_{1},\dots,Z_{k}$, where $Z_{i} = X_{[i]}$, and $[i]$ is the congruence class of $i$ modulo $n$.
  Let us denote the inductive mean of $Z_{1},\dots,Z_{k}$ with $R_{k}^{X}$.
  Similarly, we approximate $b_{n}^{Y}$ with an auxiliary sequence $W_{i}=Y_{[i]}$, and denote the inductive mean of $W_{1},\dots,W_{k}$ with $R_{k}^{Y}$.
  Setting the weights $w_{i}$ in Theorem 3.4 of \cite{lim-palfia} to $w_{i}=1/n$, for $i=1,\dots, n$, we get that $R_{k}^{X} = R_{k-1}^{X} \oplus_{1/k}Z_{k}$, and similarly $R_{k}^{Y} = R_{k-1}^{Y} \oplus_{1/k}W_{k}$. 
  The estimate of Theorem 3.4 in \cite{lim-palfia} simplifies to
  \begin{equation}\label{eq:lim-palfia-est}
    d(b_{n}^{X},R_{k}^{X})^{2} \leq \frac{n}{k} \biggl[3\Delta_{n}(X)^{2} + \frac{1}{n}\sum_{k=1}^{n}d(b_{n}^{X},X_{k})^{2}\biggr],
  \end{equation}
  where $\Delta_{n}(X) = \max_{1\leq i,j\leq n}d(X_{i},X_{j})$.
  Using the variance inequality (Lemma \ref{lemma:var-ineq}), we get
  \begin{equation*}
    \frac{1}{n}\sum_{k=1}^{n}d(b_{n}^{X},X_{k})^{2} \leq \frac{1}{n^{2}}\sum_{k=1}^{n}\sum_{j=1}^{n}d(X_{j},X_{k})^{2} \leq \Delta_{n}(X)^{2}.
  \end{equation*}
  Hence, the right-hand side of \eqref{eq:lim-palfia-est} can be estimated by
  \begin{equation*}
    d(b_{n}^{X},R_{k}^{X})^{2}\leq \frac{4 \Delta_{n}(X)^{2}n}{k}.
  \end{equation*}
  Now, setting $k=n^{2}$, the triangle inequality, and $\|x\|_{\ell^{1}}^{2}\leq 3 \|x\|_{\ell^{2}}^{2}$ for $x \in \mathbb{R}^{3}$ yields
  \begin{equation*}
  \mathbb{E}(d\bigl(b_{n}^{X},b_{n}^{Y}\bigr)^{2})\leq 3 \bigl[\mathbb{E}\bigl(d\bigl(b_{n}^{X},R_{n^{2}}^{X}\bigr)^{2}\bigr) + \mathbb{E}\bigl(d\bigl(R_{n^{2}}^{X},R_{n^{2}}^{Y}\bigr)^{2}\bigr) + \mathbb{E}\bigl(d\bigl(b_{n}^{Y},R_{n^{2}}^{Y}\bigr)^{2}\bigr)\bigr].
  \end{equation*}
  The first term can be estimated by
  \begin{equation*}
    \mathbb{E}\bigl(d\bigl(b_{n}^{X},R_{n^{2}}^{X}\bigr)^{2}\bigr) \leq \frac{4}{n}\mathbb{E}(\Delta_{n}(X)^{2}) \leq \frac{8}{n}\mathbb{E}(\max_{1\leq k \leq n}d(X_{k},a)^{2}) \to 0.
  \end{equation*}
      By a similar estimate, we get $\mathbb{E}\bigl(d\bigl(b_{n}^{X},R_{n^{2}}^{X}\bigr)^{2}\bigr) \to 0$.
      For the remaining term, we again apply \eqref{eq:d(xt,yt)} iteratively, to the effect of 
      \begin{equation*}
        \mathbb{E}(d(R_{n^{2}}^{X},R_{n^{2}}^{Y})^{2})\leq \frac{1}{n^{2}}\sum_{k=1}^{n^{2}}\mathbb{E}(d(X_{[k]},Y_{[k]})^{2}) = \frac{1}{n}\sum_{k=1}^{n}\mathbb{E}(d(X_{k},Y_{k})^{2}) \to 0,
      \end{equation*}
  and hence, $b_{n}^{Y} \to \mu$ in $\mathcal{L}^{2}(H)$.
  It is worth pointing out that a similar argument has been used to show an $L^1$ contraction property of $b_n$ in Proposition 3.1 of \cite{brunel}.

\end{proof}

\begin{proposition}\label{prob:huber-L2}
  Let $X_{1},\dots,X_{n}$, and $Y_{1},\dots,Y_{n}$ be two sequences of random variables in $\mathcal{L}^{2}(H)$. %, such that $\{X_{1},\dots,X_{n},Y_{1},\dots,Y_{n}\}$ is independent. 
  Let $f$ be uniformly drawn from
  \begin{equation*}
    \mathcal{F}=\{f:\{1,\dots,n\} \to \{0,1\}\mid |f^{-1}(1)| = r_{n}\}
  \end{equation*}
  for some rate of contamination $r_{n}\leq n$, independent of $\{X_{1},\dots,X_{n},Y_{1},\dots,Y_{n}\}$, and let
  \begin{equation*}
    Z_{k} = \begin{cases}
      X_{k} & \text{if } f(k) = 0 \\
      Y_{k} & \text{if } f(k) = 1
    \end{cases}.
  \end{equation*}
If the inductive mean of $X_{1},\dots,X_{n}$ converges to some $\mu \in H$, (e.g. if the assumptions of Theorem \ref{thm:slln-main} are satisfied) then $S_{n}^{Z}$ -- the inductive mean of $Z_{1},\dots,Z_{n}$ -- converges to $\mu$ in $\mathcal{L}^{2}(H)$, if
\begin{equation}\label{eq:huber-cond}
  \frac{r_{n}}{n^{2}}\sum_{k=1}^{n}\mathbb{E}(d(X_{k},Y_{k})^{2})\to 0.
\end{equation}
If the Fr\'echet mean of $X_1,\dots,X_n$ converges to $\mu$ in $L^2$, Assumption \eqref{eq:huber-cond} is satisfied, and both $\mathbb{E}(\max_{1\leq k\leq n} d(X_k,a)^2)$, and $\mathbb{E}(\max_{1\leq k\leq n}d(Y_k,a)^2)$ are $o(n)$ for some (and hence all) $a\in H$, then the Fr\'echet mean of $Z_1\dots,Z_n$ converges to $\mu$ in $L^2$ as well.
\end{proposition}
\begin{proof}
  Both claims follow immediately from Theorem \ref{thm:slln-main}, and Proposition \ref{prop:stab-of-two-seq}. 
  Simply note that 
  \begin{equation*}
    \begin{split}
      \mathbb{E}(d(Z_{i},X_{i})^{2}) &= \mathbb{E}(d(X_{i},X_{i})^{2} I_{\{f(i) = 0\}}) + \mathbb{E}(d(Y_{i},X_{i})^{2}I_{\{f(i)=1\}}) = \mathbb{E}(d(Y_{i},X_{i})^{2})\mathbb{P}(f(i) = 1) \\
                                     &= \mathbb{E}(d(Y_{i},X_{i})^{2})\frac{r_{n}}{n}.
    \end{split}
  \end{equation*}
\end{proof}

\begin{proposition}\label{prop:huber-as}
  Let $X_{1},\dots,X_{n}$ be independent with expectation $\mu_{k}$, and let $Y_{1},\dots,Y_{r_{n}}$ be independent with mean $\nu_{k}$, for some contamination level $r_{n}\leq n$, such that $\{X_{1},\dots,X_{n},Y_{1},\dots,Y_{r_{n}}\}$ is independent.
Let $Z_{1},\dots,Z_{n}$ be generated from $X_{1},\dots,X_{n}$ by replacing $r_{n}$ of the $X_{i}$'s with one $Y_{j}$'s (i.e., every $Y_j$ is used exactly once), and let $S_{n}$ be the inductive mean of $Z_{1},\dots,Z_{n}$.
If there is a $\mu\in H$ such that for some $p>1/2$, $C>0$, $0\leq q < \min\{1/4, p-1/2\}$, and $z\in H$ we have
\begin{enumerate}
\item $\frac{1}{n}\sum_{k=1}^{n}d(\mu_{k},\mu) = O(n^{-p})$,
\item $\frac{1}{n}\sum_{k=1}^{r_{n}}d(\nu_{k},\mu) = O(n^{-p})$,
\item $\mathbb{P}(\max_{1\leq k \leq n}d(X_{k},z) \leq Cn^{q}) = 1$, and
\item $\mathbb{P}(\max_{1\leq k \leq r_{n}}d(Y_{k},z) \leq Cn^{q}) = 1$,
\end{enumerate}
then $S_{n}\to \mu$ almost surely.
\end{proposition}

\begin{proof}
  Let $\kappa_{k}$ be the mean of $Z_{k}$. 
  It is enough to check, that the assumptions of Corollary \ref{cor:slln} are fulfilled. 
  Let us denote $I_{X}$ the set of indices for which $Z_{k} = X_{k}$.
  Note that 
  \begin{equation*}
    \frac{1}{n}\sum_{k=1}^{n}d(\kappa_{k},\mu) = \frac{1}{n}\sum_{k\in I_{X}}d(\mu_{k},\mu) + \frac{1}{n}\sum_{k=1}^{r_{n}}d(\nu_{k},\mu) \leq \frac{1}{n}\sum_{k=1}^{n}d(\mu_{k},\mu) + \frac{1}{n}\sum_{k=1}^{r_{n}}d(\nu_{k},\mu) = O(n^{-p}),
  \end{equation*}
  for some $p>1/2$.
  Let us write $M_n(X) = \max_{1\leq k \leq n}d(X_{k},z)$, $M_{n}(Y) = \max_{1\leq k \leq r_{n}}d(Y_{k},z)$, and $M_{n}(Z) = \max_{1\leq k \leq n}d(Z_{k},z)$.
  We want to show that $\mathbb{P}(M_{n}(Z) \leq C'n^{q}) = 1$ for some $C'>0$. 
  Note that 
  \begin{equation*}
    \mathbb{P}(M_{n}(Z) \leq C'n^{q}) \geq \mathbb{P}(\max\{M_{n}(X), M_{n}(Y)\} \leq C'n^{q}).
  \end{equation*}
  On the other hand, 
  \begin{equation*}
    \begin{split}
      \mathbb{P}(\max\{M_{n}(X), M_{n}(Y)\} > 2Cn^{q}) &\leq \mathbb{P}(M_{n}(X) + M_{n}(Y) > 2Cn^{q}) \\
                                                       &\leq \mathbb{P}(M_{n}(X) > Cn^{q}) + \mathbb{P}(M_{n}(Y) > Cn^{q}) = 0, 
    \end{split}
  \end{equation*}
  i.e., we can choose $C' = 2C$.
\end{proof}

\section{Means in Hadamard Spaces}\label{sec:means}
In this section, we discuss the stability of various means introduced in the literature. 
All of them allow for some version of Ces\`aro's Lemma or, if weighted versions are available, the Toeplitz Lemma. 
In other words, if the underlying sequence of points converges, means in Hadamard spaces are at least as well-behaved as means in linear spaces. 
At the end of this section, we briefly discuss the computational tractability of various approximation schemes for the Fr\'echet mean, and introduce a stochastic resampling method.
Using this resampling method, we analyze the stability of the inductive mean with respect to data loss, permutation, and noise.

Es-Sahib and Heinich \cite{es-sahib-heinich} attempted an axiomatic approach.  
On a locally compact Hadamard space $(H,d)$, one can recursively define a unique map $\beta_{n}: H^{n} \to H$ satisfying the following three axioms:
\begin{enumerate}
  \item $\beta_{n}(x,\dots,x) = x$ for all $x\in H$,
  \item $d(\beta_{n}(x_{1},\dots,x_{n}),\beta_{n}(y_{1},\dots,y_{n})) \leq \frac{1}{n}\sum_{k=1}^{n}d(x_{k},y_{k})$ for all $x_{i},y_{j}\in H$, and
  \item $\beta_{n}(x_{1},\dots,x_{n}) = \beta_{n}(\hat{x}_{1},\dots,\hat{x}_{n})$, where $\hat{x}_{i} = \beta_{n-1}(x_{1},\dots,x_{i-1},x_{i+1},\dots,x_{n})$.
\end{enumerate}
This map is symmetric and satisfies $d(z,\beta_{n}(x_{1},\dots,x_{n})) \leq \frac{1}{n}\sum_{k=1}^{n}d(z,x_{k})$.
A strong law of large numbers is available for such maps \cite{es-sahib-heinich}.
To the best of our knowledge, this was the first law of large numbers specifically for Hadamard spaces (for early results in general metric spaces see \cite{ziezold1977,slln-centroid-cpt-metric-space}). 
This mean is in general different from the barycenter $b_{n}$ and the inductive mean $S_{n}$.
In particular $\beta_{n}$ and $b_{n}$ are $L^{1}$-contracting \cite{sturm}.
On the other hand, iterating the geodesic comparison lemma (Corollary 2.5 in \cite{sturm}), we get that $S_{n}$ is $L^{2}$-contracting on finitely-supported probability measures. 
Since $d(\cdot,\cdot)^{2}$ is convex on $H^{2}$ by the NPC inequality \eqref{eq:npc-ineq}, the proof of Theorem 6.3 in \cite{sturm} can be adapted to show that $b_{n}$ is also $L^{2}$-contracting (this is also established as a by-product of the proof of Proposition \ref{prob:huber-L2}).
In \cite{navas-ergodic} Navas generalized the construction of Es-Sahib and Heinich to potentially non-locally compact spaces.
His goal was to use a well-behaved mean to establish an $L^{1}$ ergodic theorem. 
The notion of a \textit{convex mean} \cite{sturm} of a random variable (or probability distribution) constitutes yet another approach.
Given $X \in \mathcal{L}^{1}(H)$ we say that $z \in H$ is a convex mean of $X$, if for all convex, Lipschitz continuous functions $\varphi: H\to \mathbb{R}$ we have $\varphi(z) \leq \mathbb{E}(\varphi(X))$.

Hansen \cite{Hansen} has introduced a mean specifically on the space of positive definite symmetric matrices. 
His goal was to construct a computationally more tractable mean than the Fr\'echet mean $b_{n}$. 
In this sense, the mean of Hansen is similar to $S_{n}$. 
Hansen's mean is based on the following observation: in Euclidean spaces, the mean of $n$ points can also be rewritten as 
\begin{equation*}
  \frac{1}{n}\sum_{k=1}^{n}x_{k} = \frac{\bigl(\frac{n-1}{n}x_{1} + \frac{1}{n}x_{n}\bigr) + \cdots + \bigl(\frac{n-1}{n}x_{n-1} + \frac{1}{n}x_{n}\bigr)}{n-1},
\end{equation*}
i.e., we have a convex combination of convex combinations that can be computed recursively.
Translating this idea into the language of Hadamard spaces (see \cite{kim}), we may inductively define $H_{n} = H_{n}(x_{1},\dots,x_{n})$ by
\begin{equation*}
  H_{1} = x_{1}, \quad \text{and} \quad H_{n} = H_{n-1}(x_{1}\oplus_{\frac{1}{n}}x_{n},\dots, x_{n-1}\oplus_{\frac{1}{n}}x_{n}).
\end{equation*}
In the Euclidean case $S_{n}$, $H_{n}$, and $b_{n}$ coincide, but in a Hadamard space they are known to differ (\cite{kim}, Example 5.3). 
Now the obvious question arises: How close are these means to each other in general?
It turns out, that if the underlying sequence of points converges, the sequence of means converges to the same point for all means introduced in this section.
In other words, Ces\`aro's lemma is a universal property of means in spaces of non-positive curvature.
\begin{lemma}[Stability of Ces\`aro means]\label{lemma:cesaro}
  Let $(x_{n})_{n\geq 1}$ be a sequence of points in a Hadamard space converging to $x$. 
  Furthermore, let $b_{n}$ be their barycenter, $S_{n}$ be the inductive mean, $\beta_{n} = \beta_{n}(x_{1},\dots,x_{n})$ be the mean introduced by Navas, and let $c_{n}$ be any convex mean of $X_{n}$, where $X_{n}$ is drawn uniformly from $\{x_{1},\dots,x_{n}\}$. Then $b_{n},S_{n},\beta_{n},c_{n} \to x$.
\end{lemma}
\begin{proof}
  Starting with $b_{n}$, we observe that simply neglecting the negative part in the variance inequality \eqref{eq:var} yields
  \begin{equation*}
    d(x,b_{n})^{2} \leq \frac{1}{n}\sum_{k=1}^{n}d(x,x_{k})^{2}.
  \end{equation*}
  By the classical Ces\`aro lemma, this implies $b_{n} \to x$. 
  Using Inequality \eqref{eq:sn-est} gives
  \begin{equation*}
    d(x,S_{n})^{2} \leq \frac{1}{n}\sum_{k=1}^{n}d(x,x_{k})^{2} \to 0.
  \end{equation*}
Moving on to the mean of Navas, we note that the first two axioms imply
\begin{equation*}
  d(x,\beta_{n}) = d(\beta_{n}(x,\dots,x),\beta_{n}(x_{1},\dots,x_{n})) \leq \frac{1}{n}\sum_{k=1}^{n} d(x,x_{k}),
\end{equation*}
which again goes to $0$ by the classical Ces\`aro Lemma. For the convex mean $c_{n}$, we note that $z\mapsto d(x,z)^{2}$ is convex and Lipschitz continuous. Hence,
\begin{equation*}
  d(x,c_{n})^{2} \leq \frac{1}{n}\sum_{k=1}^{n}d(x,x_{k})^{2} \to 0.
\end{equation*}
\end{proof}

Choi and Ji \cite{choi-toeplitz} proved a version of the classical Toeplitz Lemma (p. 250, \cite{loeve}), for weighted versions of $S_{n}$ and Hansen's mean.
Lemma \ref{lemma:toeplitz} below mimics their result for weighted versions of $b_{n}$. 
This suggests, that taking averages on Hadamard spaces is about as stable as in the Euclidean case if the underlying sequence of points converges. 
\begin{lemma}\label{lemma:toeplitz}
Let $p^{(n)} = (p^{(n)}_{1}, \dots, p_{n}^{(n)})$ be a sequence of probability measures on $n$ points, such that $\lim_{n\to\infty}p^{(n)}_{k} = 0$ for all $k\geq 1$, and let $(x_{n})_{n\geq 1}$ be a sequence of points in a Hadamard space $(H,d)$, such that $x_{n} \to x$. 
If $X_{n}$ is drawn from $\{x_{1},\dots,x_{n}\}$ randomly with $\mathbb{P}(X_{n}=x_{k}) = p^{(n)}_{k}$, then $\mathbb{E}(X_{n}) \to x$.  
\end{lemma}
\begin{proof}
  Ignoring the second term in the variance inequality \eqref{eq:var} yields
  \begin{equation*}
    d(\mathbb{E}(X_{n}),x)^{2} \leq \sum_{k=1}^{n}p^{(n)}_{k}d(x_{k},x)^{2},
  \end{equation*}
  which goes to zero, by the classical Toeplitz lemma (p. 250, \cite{loeve}).
\end{proof}

Lemma \ref{lemma:cesaro} and Lemma \ref{lemma:toeplitz} show that, if the underlying sequence of points converges, means on Hadamard spaces are as well behaved as in linear spaces. 
However, the computation of these means is much more challenging than in the linear case.
Various attempts have been made to approximate theoretically well-behaved means by computationally tractable ones.
For the computation of the barycenter $b_{n}$ of $n$ points Lim and P\'alfia \cite{lim-palfia} proposed the following scheme. 
If $y_{k} = x_{[k]}$, where $[k]$ is the residue of $k$ modulo $n$, we set 
\begin{equation*}
  LP_{n}^{(1)} = y_{1}, \quad \text{and} \quad LP_{n}^{(k)} = LP_{n}^{(k-1)} \oplus_{\frac{1}{k}} y_{k}.
\end{equation*}
In other words, Lim and P\'alfia \cite{lim-palfia} compute the inductive mean of the $n$-periodic sequence
\begin{equation*}
  x_{1},\dots,x_{n},x_{1},\dots,x_{n},x_{1},\dots.
\end{equation*}
For fixed $n$, they show that $LP_{n}^{(k)} \to b_{n}$ as $k \to \infty$. 
If the underlying set of points $\{x_{n}\}_{n\in \mathbb{N}}$ is bounded, Theorem 3.4 in \cite{lim-palfia} implies that for any $n$ and $k$
\begin{equation*}
  d(LP_{n}^{(k)}, b_{n}) \leq 2\Delta\sqrt{\frac{n}{k}},
\end{equation*}
where $\Delta$ is the diameter of the points $\{x_{n}\}_{n\in \mathbb{N}}$.
If we naively use $LP_{n}^{(k)}$ to approximate $\mu$, we end up with
\begin{equation}\label{eq:lim-palfia-rate}
  d(LP_{n}^{(k)},\mu) \leq d(LP_{n}^{(k)},b_{n}) + d(b_{n},\mu) \leq 2 \Delta \sqrt{\frac{n}{k}} + d(b_{n},\mu).
\end{equation}
This has two drawbacks. 
First, we need to compute $nf(n)$ geodesics to arrive at an error $d(LP_{n}^{(nf(n))},b_{n}) = O(f(n)^{-1/2})$ for any $f:\mathbb{N}\to \mathbb{N}$, such that $\lim_{n\to \infty} f(n) = \infty$.
This can be computationally expensive.  
Additionally, even if we can efficiently compute $LP_{n}^{(nf(n))}$, we may waste a lot of computational effort, as the overall error still depends on the term $d(b_{n},\mu)$, of which we do not know the rate of convergence in practice.
Hence, we may spend disproportional effort minimizing the first term of the right-hand side of \eqref{eq:lim-palfia-rate}, without minimizing the overall error. 
Yet another algorithm to approximate $b_n$ is the proximal point algorithm, which is popular for optimization problems in Hilbert spaces, and can be extended to the Hadamard space setting (cf. chapter 6.3, \cite{bacakbook}). Starting at an arbitrary point $y_0 \in H$, the $m$-th step of the algorithm is given by 
\begin{equation}
  y_m \coloneqq  \underset{z \in H}{\mathrm{arg min}}\bigg\{\tfrac{1}{n}\sum_{i=1}^n d(x_i,z)^2+\tfrac{1}{2\lambda_m}d(y_{m-1},z)^2\bigg\}, \quad \lambda_m>0,\,\,\text{such that}\,\, \sum_{k=1}^{\infty}\lambda_k=\infty .
\end{equation} 
The sequence $(y_m)_{m\geq1}$ converges \textit{weakly} (in the sense of Definition 3.1.1 in \cite{bacakbook}) to the unique minimizer of $n^{-1}\sum_{i=1}^n d(x_i,z)^2$ (Theorem 6.3.1, \cite{bacakbook}).
Finding the minimizer in each step can be computationally cumbersome, and does not provide any information about $d(b_n,\mu)$.

Ideally, we would like to use the inductive mean $S_n$ as a proxy for $b_n$.
Theorem \ref{thm:slln-main} and Theorem \ref{thm:limit-bn}, together with Propositions \ref{prob:huber-L2} and \ref{prop:huber-as} provide some justification for this. 
However, since the inductive mean is not symmetric, neither is the approximation scheme of Lim and P\'alfia \cite{lim-palfia}. 
This lack of symmetry has to be taken into account.
In particular, if $S_n$ is meant to be a good proxy for $b_n$, we would expect $S_n$ to be \textit{asymptotically symmetric}, in some sense. 
To address the issues of data loss, non-symmetry, and contamination at the same time, we propose the following stochastic approximation scheme.
Let $Y_{k}$ be drawn independently and uniformly from the set $\{x_{1},\dots,x_{k}\}$, and define
\begin{equation*}
  M_{1} = Y_{1}, \quad \text{and} \quad M_{n} = M_{n-1} \oplus_{\frac{1}{n}} Y_{n}.
\end{equation*}
This approximation scheme is justified using our heteroscedastic law of large numbers (Corollary \ref{cor:slln}). 
Note that the $Y_i$'s are not identically distributed.
However, their means corresponds to the Fr\'echet means of the first $i$ points. 
\begin{theorem}\label{thm:approx-scheme}
  Let $\{x_{n}\}_{n\in \mathbb{N}}$ be a set of points in a Hadamard space $(H,d)$, and let $b_{n}$ be the barycenter of $x_{1},\dots,x_{n}$. If there is a $\mu\in H$, such that 
  \begin{equation*}
    \frac{1}{n}\sum_{k=1}^{n}d(b_{k},\mu) = O(n^{-p})
  \end{equation*}
for some $p>1/2$, and if there is a $z \in H$, such that $x_{1},\dots,x_{n} \in B_{Cn^{q}}(z)$ for some $C>0$ and $0 \leq q < \min\{\frac{1}{4}, p-\frac{1}{2}\}$, then $M_{n} \to \mu$ almost surely.
\end{theorem}
\begin{proof}
  Since $\mathbb{E}(Y_{n}) = b_{n}$, and the sequence $(Y_{n})_{n\geq 1}$ is supported on $B_{Cn^{q}}(z)$ almost surely, Corollary \ref{cor:slln} applies.
\end{proof}
\begin{remark}
Note that we do not require the sequence of points to converge. 
In fact, not even the sequence of their Fr\'echet means does have to converge. 
The conditions imposed on the sequence are much weaker than the ones in Lemma \ref{lemma:cesaro}.
This suggests, that $b_{n}$ does not necessarily has to converge, for it to be close to $S_{n}$.
\end{remark}
\begin{remark}
    The idea of this scheme is similar to the bootstrap \cite{efron79}.
Indeed, one may ask why not to randomly sample $Y_{i}$ from the set $\{x_{1},\dots,x_{n}\}$ instead of $\{x_1,\dots,x_i\}$, and perform an \textit{n-out-of-n} bootstrap.
One of the advantages of the inductive mean over the Fr\'echet mean is the fact, that it can be updated in an online fashion. 
Sampling $Y_{i}$ from $\{x_{1},\dots,x_{i}\}$ ensures that our scheme  inherits this online updateability. 
Hence, our procedure is effective in a setting where data only becomes gradually available, and thus maintains one of the core advantages of the inductive mean. 
\end{remark}

Let us end this section by providing an overview of different means, and algorithms to compute them.
\begin{center}
    \begin{tabular}{C{2cm}C{8cm}C{5cm}}
\toprule
   Mean  &  Theory &  Computation \\
   \midrule
   Fr\'echet mean & $O(n^{-1/2})$ $L^2$-concentration requires curvature bounds \cite{le-gouic-fast-convergence}, subexponential distributions \cite{escande2023}, or entropy conditions (with either bounded spaces \cite{entropy-ahidarcoutrix} or moment conditions \cite{schötz18-entropy}). 
   $O(n^{-1/2})$ high probability results under curvature bounds and sub-Gaussian distributions in \cite{brunel} (Gaussian concentration), or sub-exponential distributions in \cite{schötz18-entropy} (exponential concentration).   & Proximal point algorithm \cite{bacakbook} (weak convergence), or approximation scheme of Lim and P\'alfia\cite{lim-palfia}, with runtime $O(nf(n)^2)$ for $O(f(n)^{-1})$ error. In the case of positive-definite, symmetric matrices strongly converging iteration schemes are available \cite{bini2013, lawson-lim}.\\
   \midrule
   Inductive mean & $O(n^{-1/2})$ $L^2$-concentration, and \& SLLN \cite{sturm}, $L^1$ ergodic theorem \cite{antezana2023}. & Direct computation via recursion in $O(n)$. \\
   \midrule
   Hansen's mean & $L^2$-LLN and SLLN in \cite{choi-toeplitz}. Finite sample behavior for positive-definite, symmetric matrices in \cite{kim}. & Direct computation via recursion in $O(n^2)$.\\
   \midrule
   Axiomatic means of \cite{es-sahib-heinich,navas-ergodic} & SLLN and $L^1$ ergodic theorems \cite{es-sahib-heinich,navas-ergodic}. & Based on a limiting process. In general, explicit error bounds for the approximation (as in \cite{lim-palfia}) are not available.\\
   
\bottomrule
\end{tabular}
\end{center}

\section{Simulations}\label{sec:sim}
We consider two examples in this simulation study. First, we look at a sequence of $2\times2$ diagonal matrices
\begin{equation*}
  A_{n} = \begin{pmatrix}
    \frac{1}{10} + \frac{1}{n} & 0 \\
    0 & 10+\frac{1}{n}
  \end{pmatrix}.
\end{equation*}
This sequence is contaminated, using our contamination model from above, with a noise matrix $B = 5I$, where $I$ denotes the $2\times2$ identity matrix, i.e., a fixed percentage of the sequence $A_{n}$ is randomly replaced with the noise matrix $B$.
We consider different levels of contamination until no estimator can recover the signal.  
The limit of $A_{n}$ is denoted with $A$.
In the case of commuting matrices $A_{1},\dots, A_{n}$, the barycenter $b_{n}$ is given by $(A_{1}\dots A_{n})^{\frac{1}{n}}$ \cite{bini2013}, $b_n = S_n = LP^{(n^2)}$.
Applying this to the sequence $A_{n}$ yields
\begin{equation*}
  b_{n} = S_{n} = LP^{(n^{2})}_{n} = \begin{pmatrix}
    \prod_{j=1}^{n}\bigl(1/10 + 1/j\bigr)^{\frac{1}{n}} & 0 \\
    0 & \prod_{j=1}^{n}\bigl(10 + 1/j\bigr)^{\frac{1}{n}}
  \end{pmatrix}.
\end{equation*}
%Since the unique minimal geodesic $A\oplus_{t}B$  between two positive-definite matrices $A$ and $B$ is given by 
%\begin{equation*}
 % A \oplus_{t} B = A^{1/2}(A^{-1/2}BA^{-1/2})^{t}A^{1/2},
%\end{equation*}
%(see \cite{kim}) a direct computation yields that the barycenter $b_{n}$, the inductive mean $S_{n}$, and the approximation scheme of Lim and P\'alfia $LP^{(n^{2})}_{n}$ all agree.
%In particular, for diagonal matrices, we have
%\begin{equation*}
%  \begin{pmatrix}
%    a_{1} & 0 \\
%    0 & a_{2}
%\end{pmatrix}\oplus_{t} 
%\begin{pmatrix}
%  b_{1}& 0 \\
%  0 & b_{2}
%\end{pmatrix} = 
%\begin{pmatrix}
%  a_{1}^{t-1}b_{1}^{t}& 0 \\
%  0 & a_{2}^{t-1}b_{2}^{t}
%\end{pmatrix}.
%\end{equation*}
However, Hansen's mean $H_{n}$ may be different from any of these.
The intrinsic metric on the space of positive definite matrices is given by $d_{P}(A,B) = \|\log(B^{-1/2}AB^{-1/2})\|_{F}$.
In addition to the intrinsic metric, we compute the distance between the estimators and their target in spectral norm. 
This shows, that comparing non-linear estimators in spectral norm can be quite misleading. 
While we obviously have
\begin{equation*}
  \|A-B\|_{2}\leq \|A-B\|_{F}\leq d_{P}(A,B),
\end{equation*}
sequences close in spectral norm may significantly differ with respect to $d_{P}$ in practice. 

\begin{figure}
\centering
\begin{subfigure}[b]{0.475\textwidth}
  \centering
  \includegraphics[width=\textwidth]{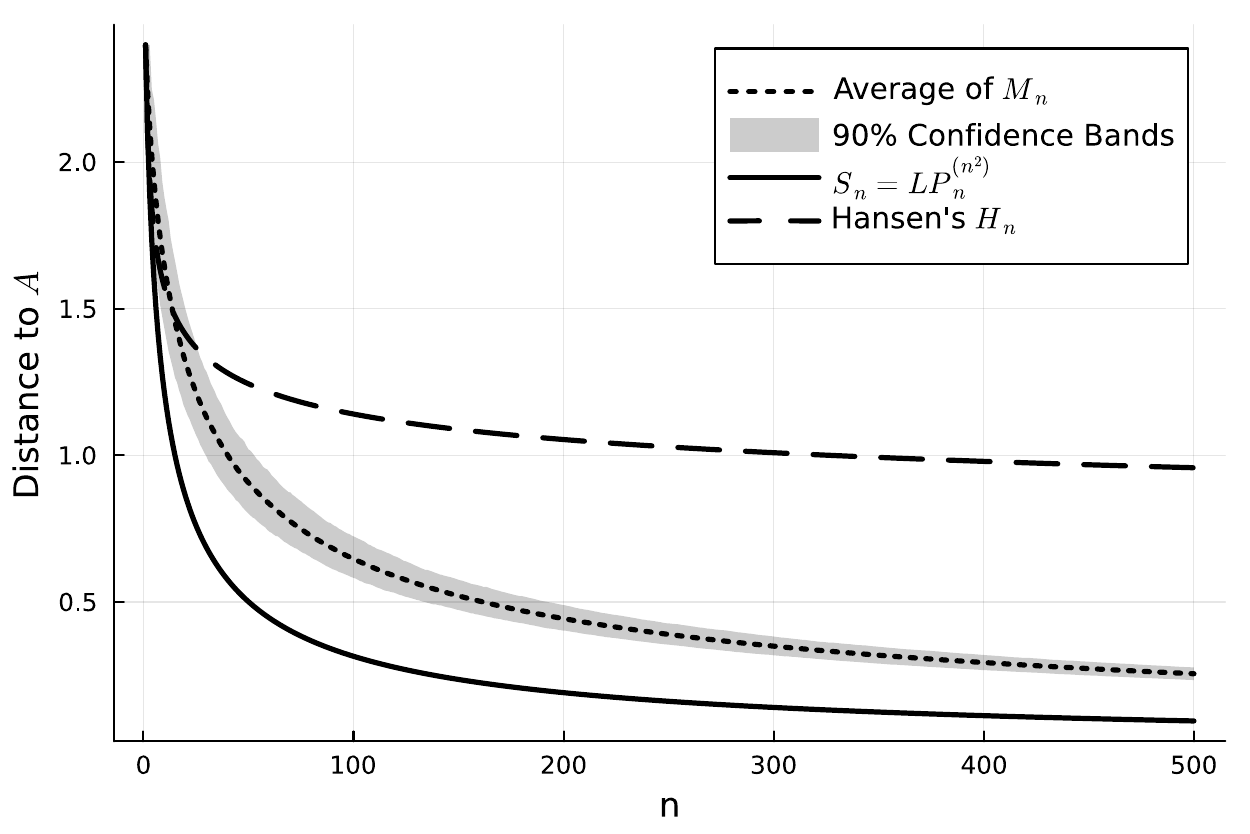}
  \caption{0\% Contamination.}
\end{subfigure}
\hfill
\begin{subfigure}[b]{0.475\textwidth}
  \centering
  \includegraphics[width=\textwidth]{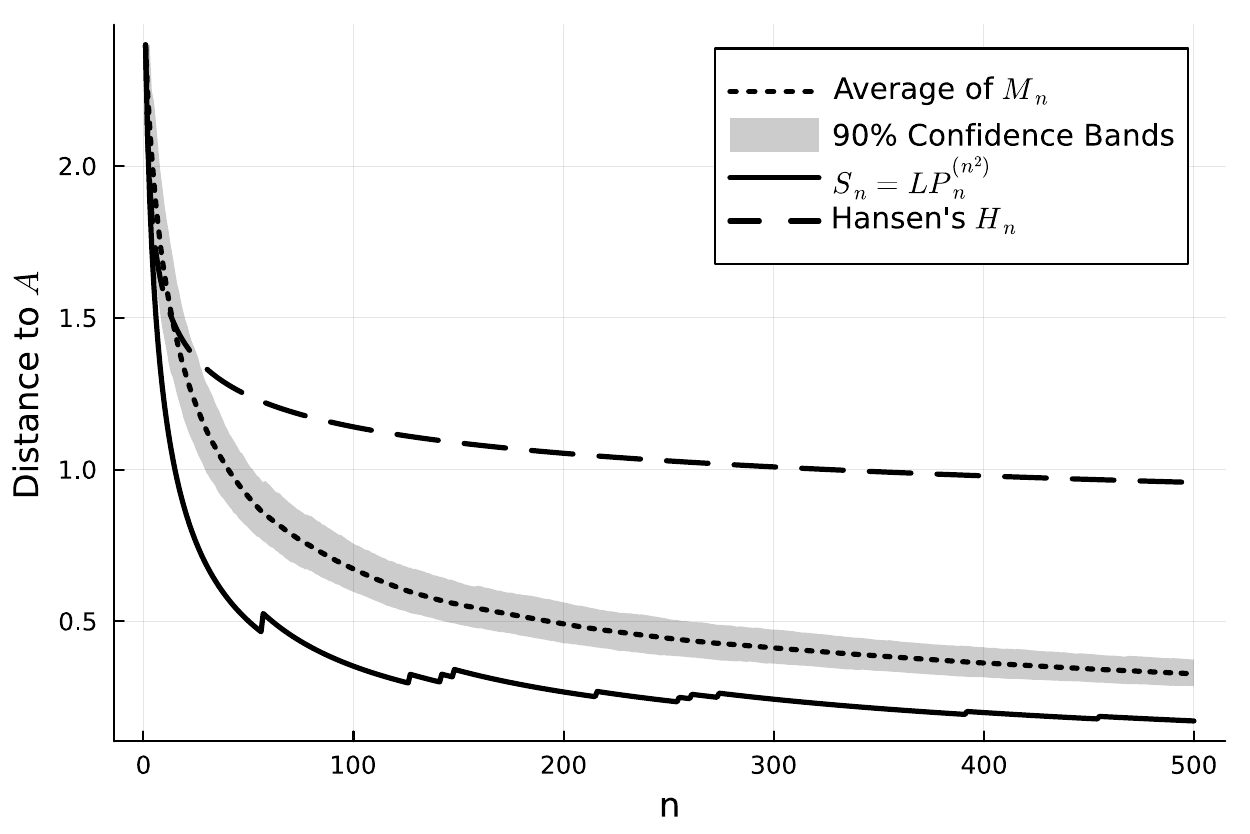}
  \caption{2\% Contamination.}
\end{subfigure}
\vskip\baselineskip
\begin{subfigure}[b]{0.475\textwidth}
  \centering
  \includegraphics[width=\textwidth]{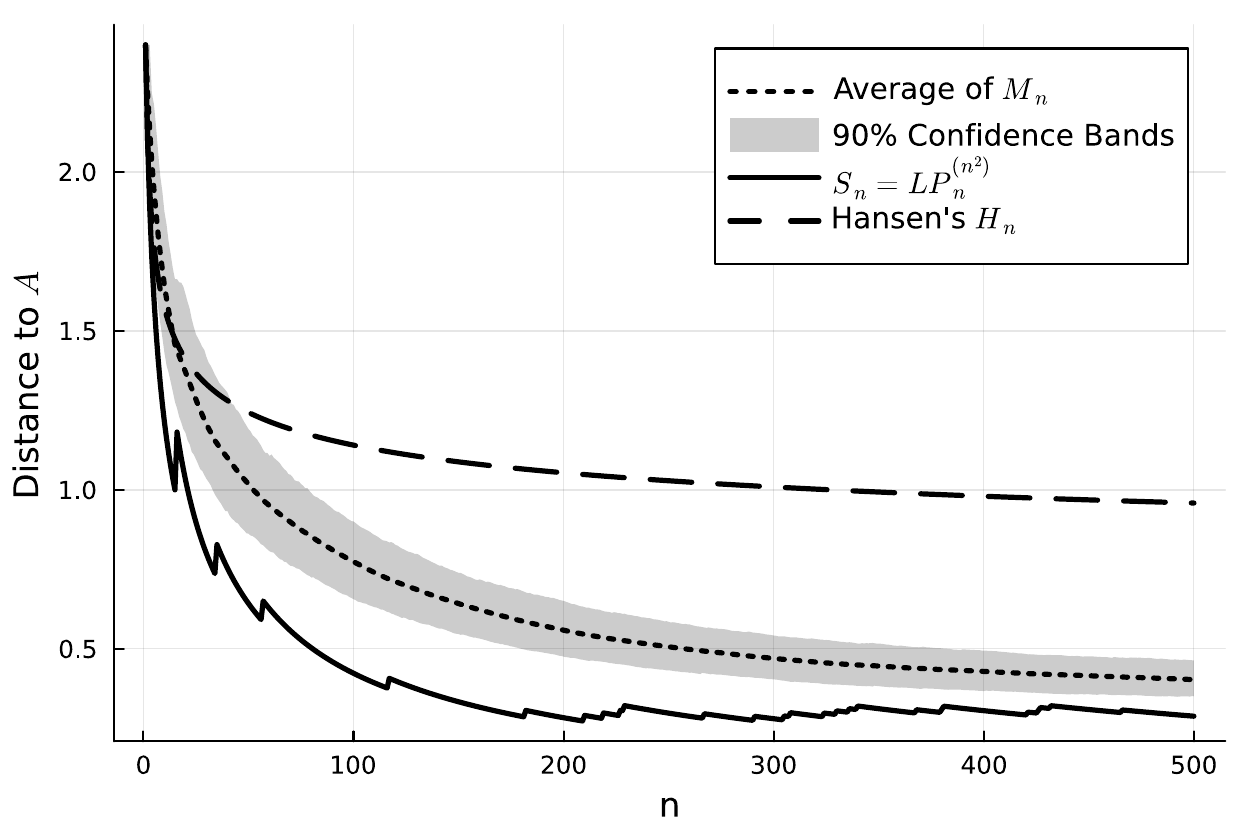}
  \caption{5\% Contamination.}
\end{subfigure}
\hfill
\begin{subfigure}[b]{0.475\textwidth}
  \centering
  \includegraphics[width=\textwidth]{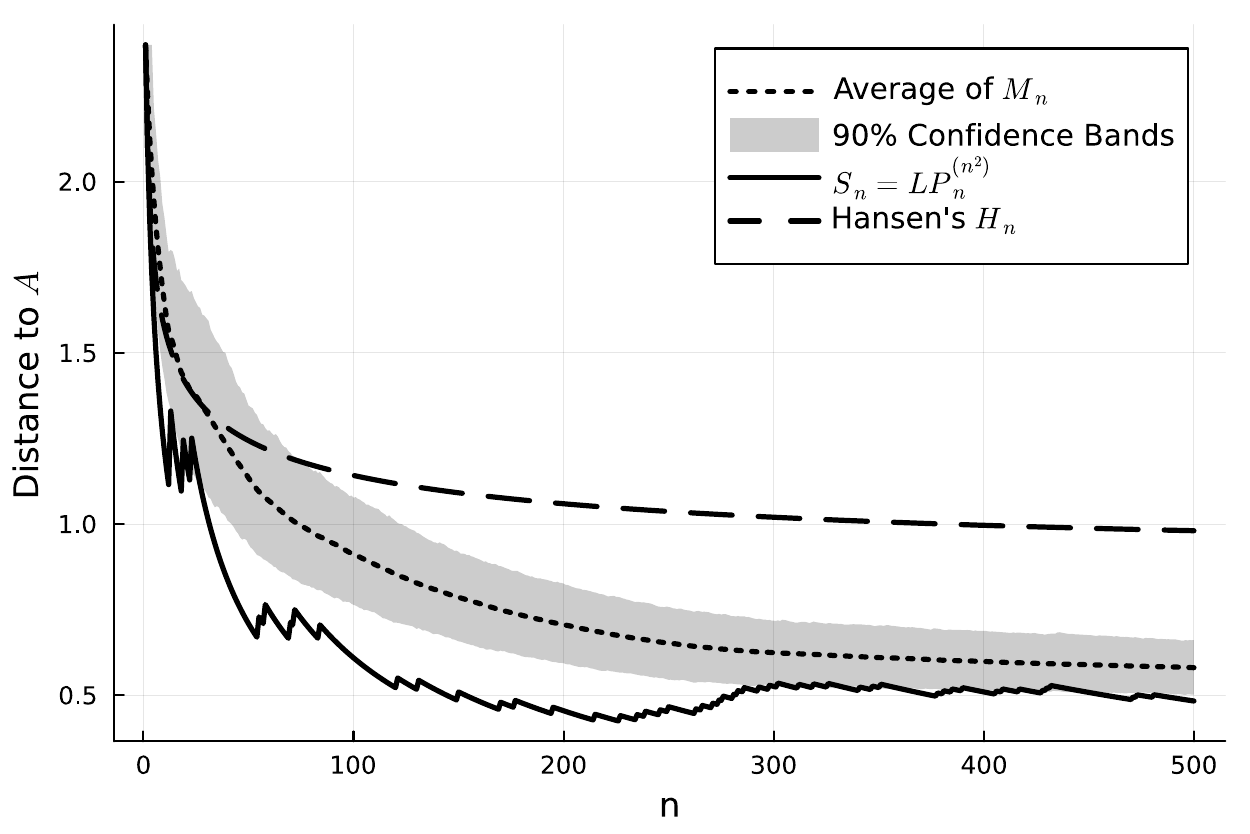}
  \caption{10\% Contamination.}
\end{subfigure}
\caption{Metric distance of $M_{n}, S_{n}$ and $H_{n}$ to $A$  at various levels of contamination.}
\label{fig:diag-metric}
\end{figure}

Comparing $S_{n},H_{n}$, and $M_{n}$ to $A$ in the metric $d_{P}$ we see that $S_{n}$ and $M_{n}$ recover the signal asymptotically, while $H_{n}$ fails to do so. 
Unsurprisingly, $S_{n}$ consistently outperforms our estimator $M_{n}$, as $M_{n}$ throws information away that $S_{n}$ can utilise.
Nonetheless, $M_{n}$ recovers $A$ asymptotically, albeit at a somewhat slower pace.
However, it is surprising, that $H_{n}$ fails to recover $A$ entirely. 
In particular, empirical observations suggest that $H_{n}$ converges to a matrix $H^{*}$ with
\begin{equation*}
  H^{*} \approx \begin{pmatrix}
    0.26 & 0 \\
    0 & 10.12
  \end{pmatrix}
\end{equation*}
under various degrees of contamination. This limit seems to be more stable with respect to contamination than $A$. However, Hansen's $H_{n}$ does not recover the \textit{signal} $A$.  
\begin{figure}
  \centering
\begin{subfigure}[b]{0.475\textwidth}
  \centering
  \includegraphics[width=\textwidth]{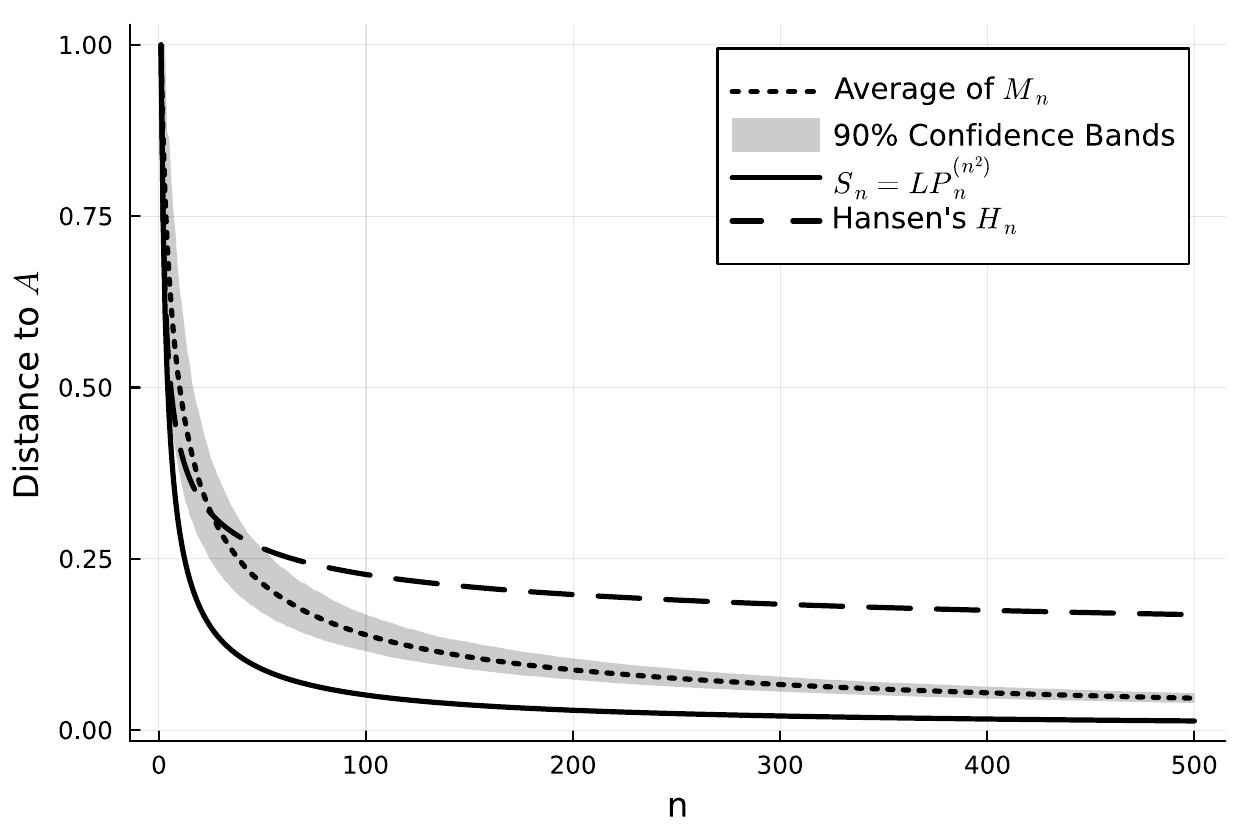}
  \caption{0\% Contamination.}
\end{subfigure}
\hfill
\begin{subfigure}[b]{0.475\textwidth}
  \centering
  \includegraphics[width=\textwidth]{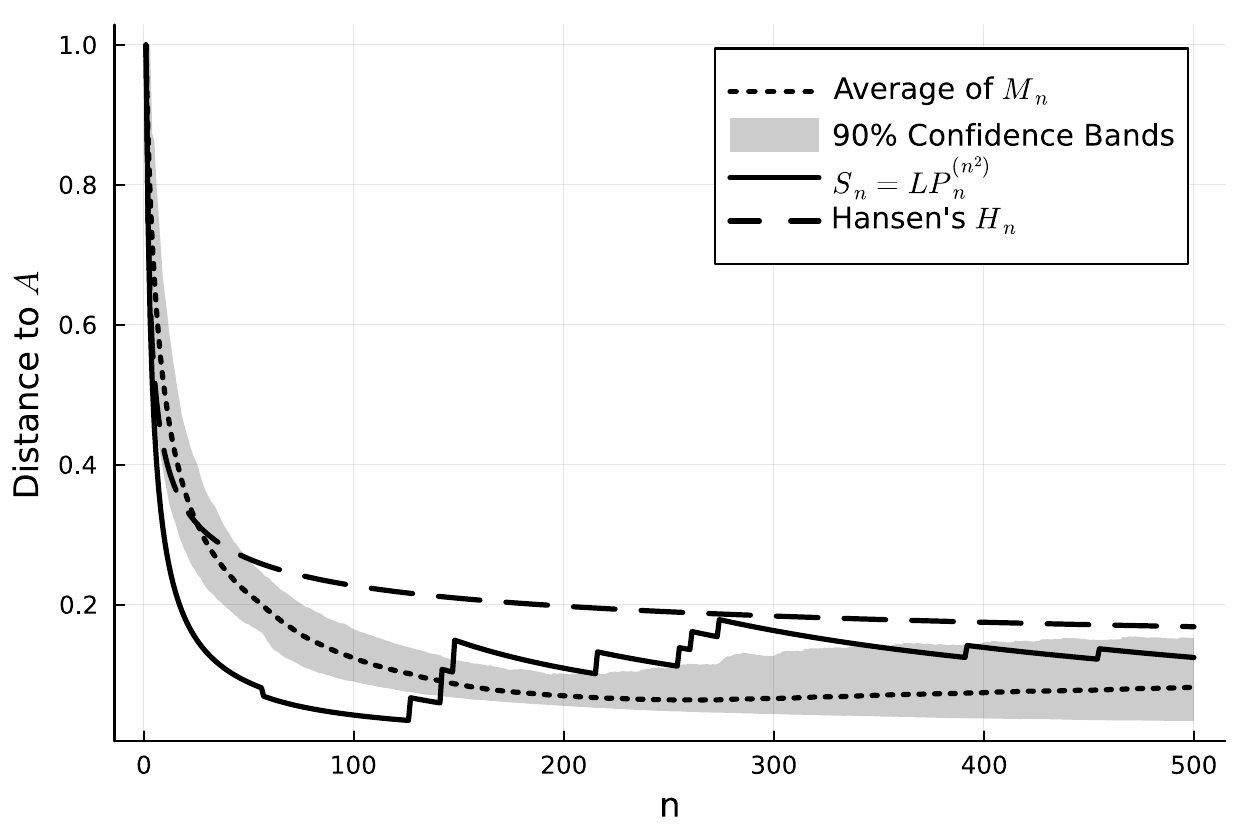}
\caption{2\% Contamination.}
\end{subfigure}
\vskip\baselineskip
\begin{subfigure}[b]{0.475\textwidth}
  \centering
  \includegraphics[width=\textwidth]{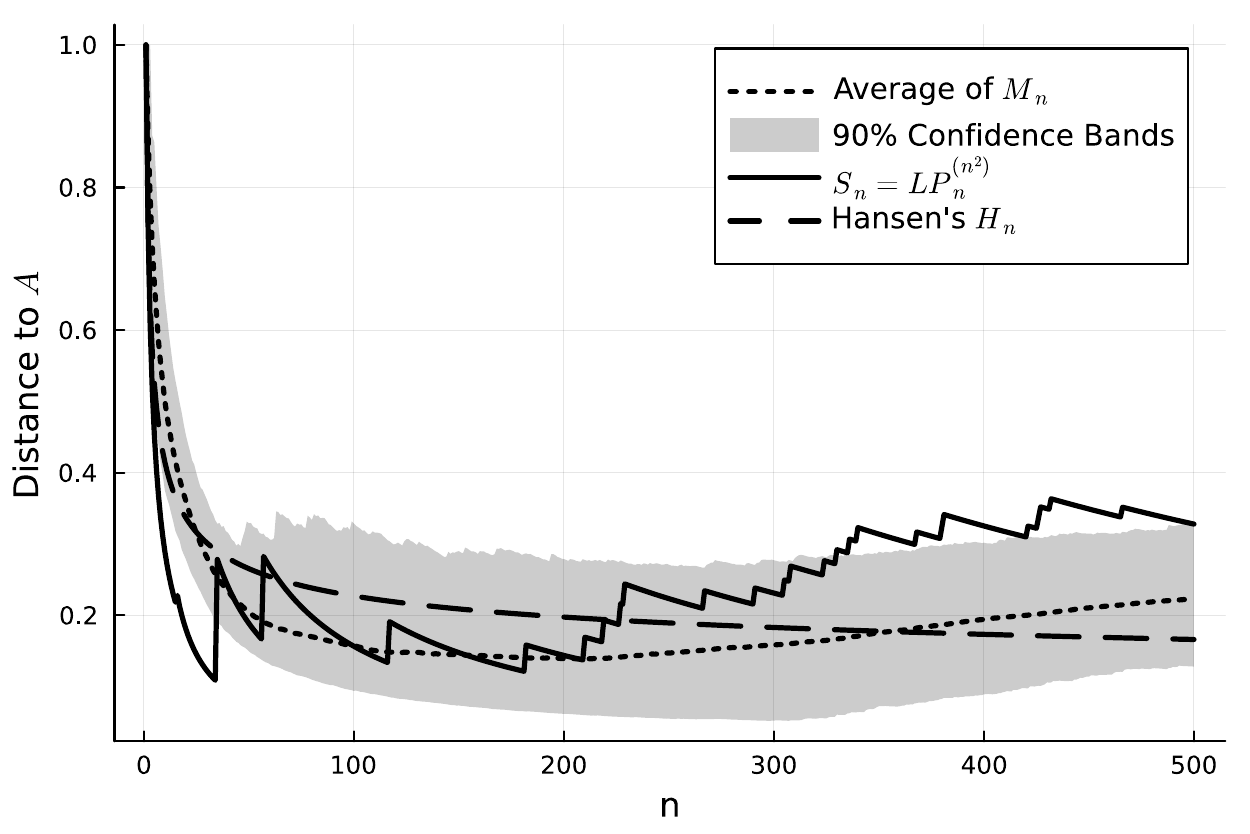}
  \caption{5\% Contamination.}
\end{subfigure}
\hfill
\begin{subfigure}[b]{0.475\textwidth}
  \centering
  \includegraphics[width=\textwidth]{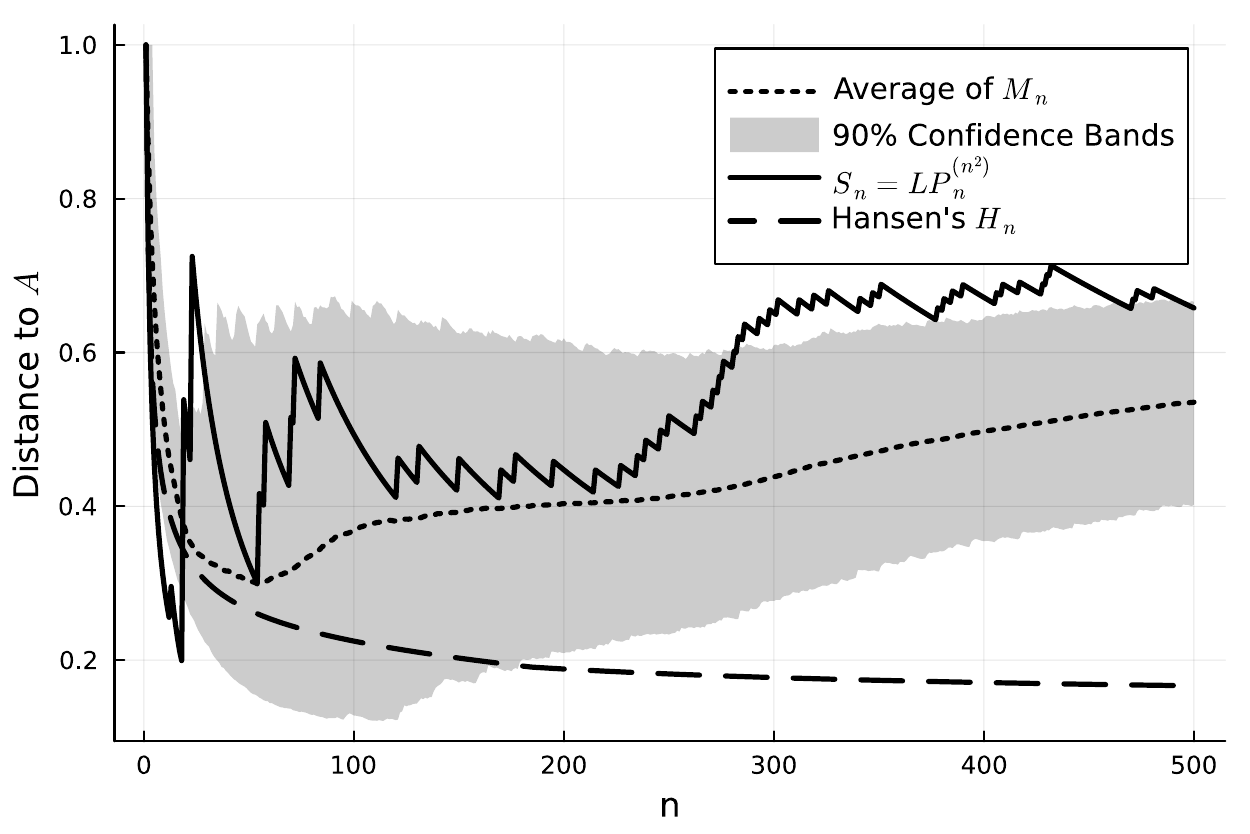}
  \caption{10\% Contamination.}
\end{subfigure}
\caption{Distance of $M_{n}, S_{n}$ and $H_{n}$ to $A$ in spectral norm at various levels of contamination.}
\label{fig:diag-norm}
\end{figure}

Measuring the distance from $S_{n}, H_{n}$, and $M_{n}$ to $A$ in spectral norm (Figure \ref{fig:diag-norm}) rather than $d_{p}$ reveals some interesting phenomena. 
All estimators seem to be more sensitive to noise when measuring the distance in spectral norm. 
Our estimator $M_{n}$ outperforms $S_{n}$ and $H_{n}$, once a certain degree of the data is contaminated, which is unexpected since we randomly throw away information. 
At about 5\% contamination, none of the above estimators can recover the signal asymptotically. 
This has practical relevance since many real-world datasets contain more than 5\% of contaminated data. 
\begin{figure}\label{fig:book}
\centering
\begin{subfigure}[b]{0.475\textwidth}
  \centering
  \includegraphics[width=\textwidth]{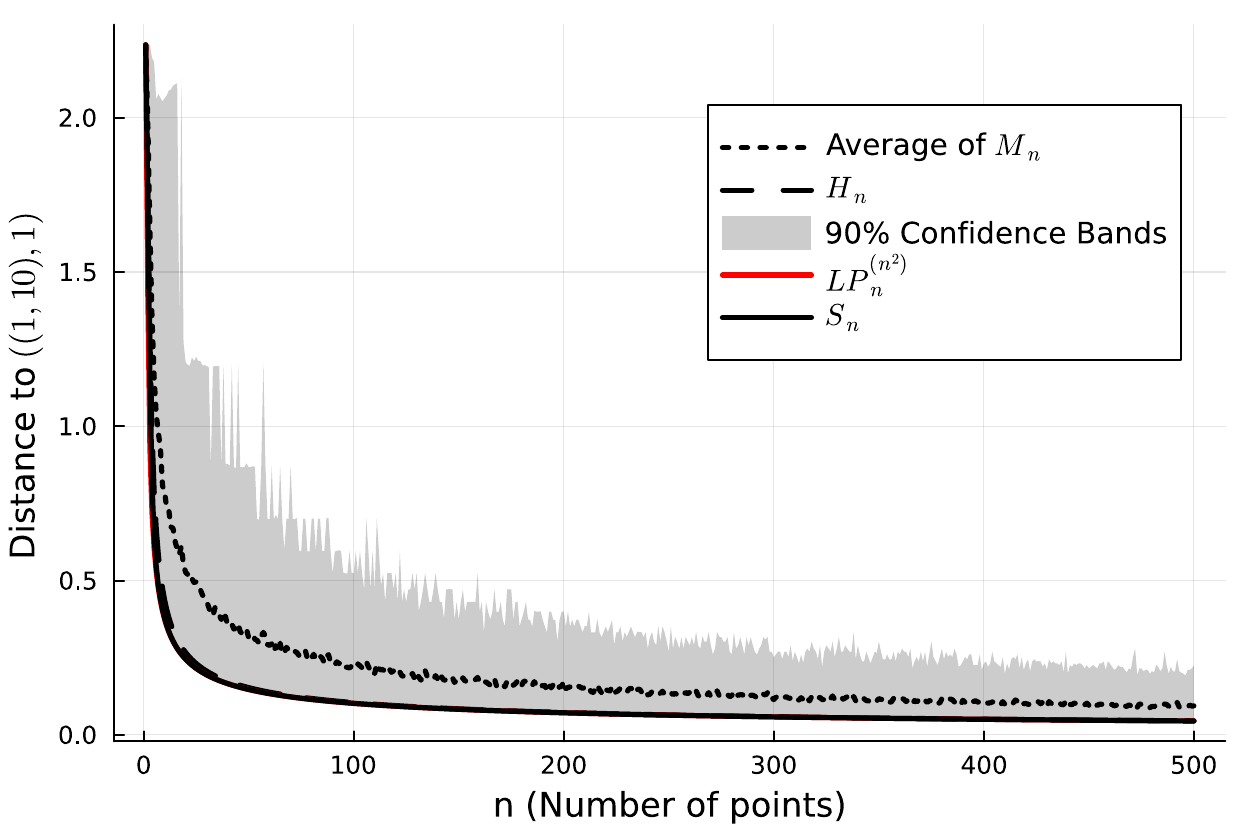}
  \caption{0\% Contamination.}
\end{subfigure}
\hfill
\begin{subfigure}[b]{0.475\textwidth}
  \centering
  \includegraphics[width=\textwidth]{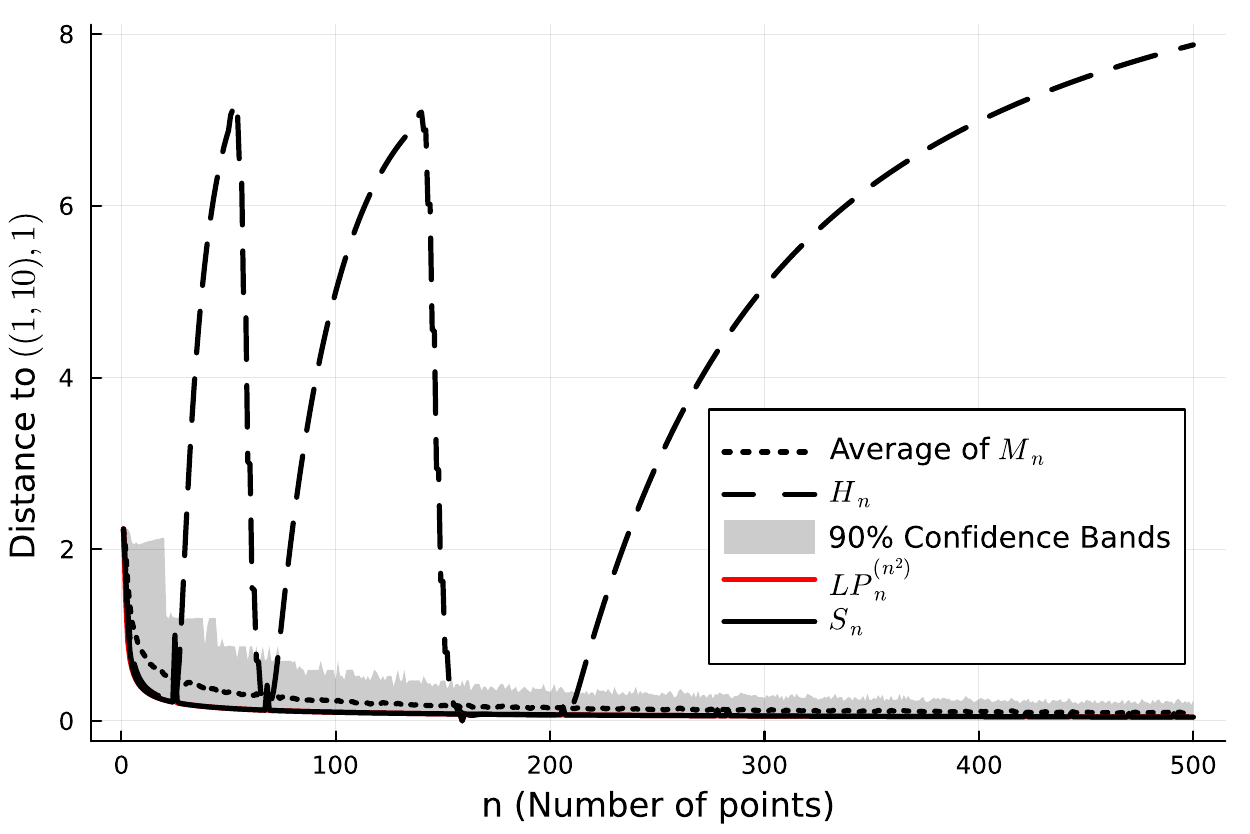}
\caption{2\% Contamination.}
\end{subfigure}
\vskip\baselineskip
\begin{subfigure}[b]{0.475\textwidth}
  \centering
  \includegraphics[width=\textwidth]{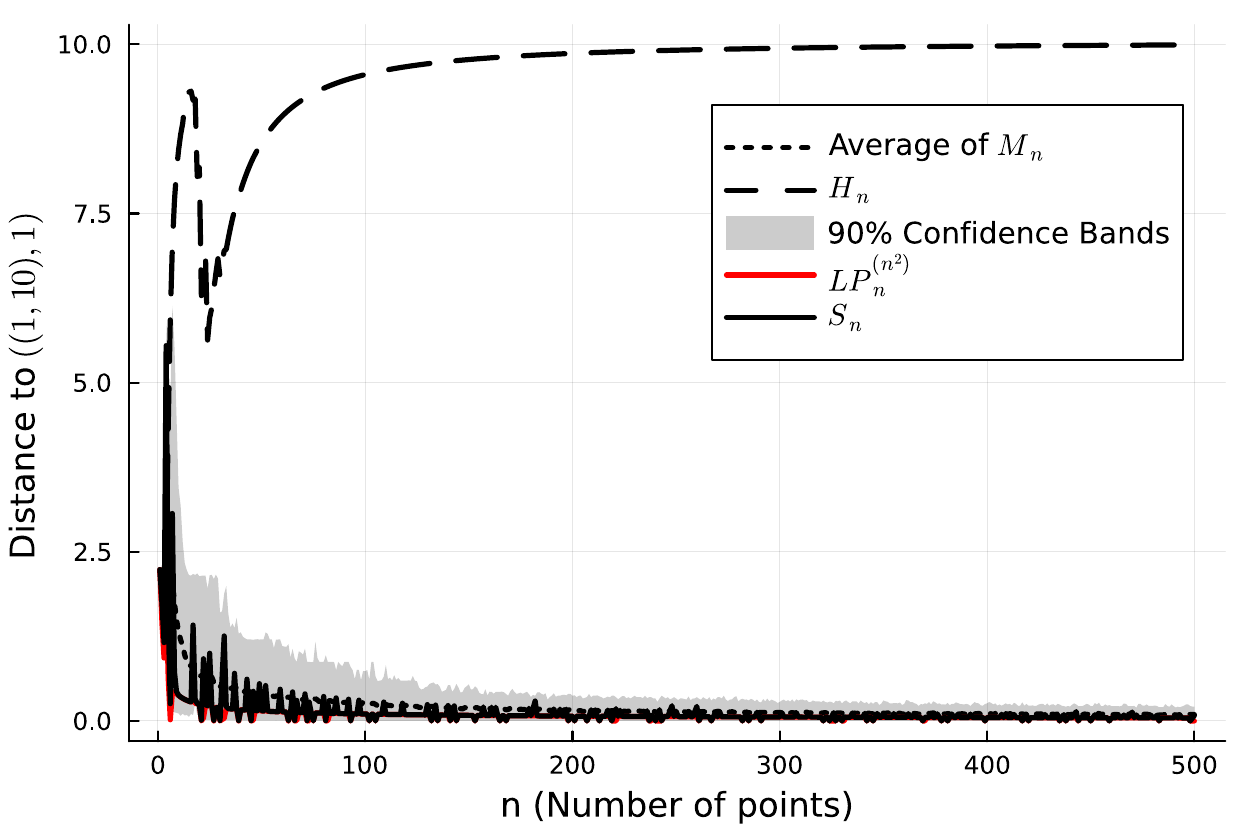}
  \caption{30\% Contamination.}
\end{subfigure}
\hfill
\begin{subfigure}[b]{0.475\textwidth}
  \centering
  \includegraphics[width=\textwidth]{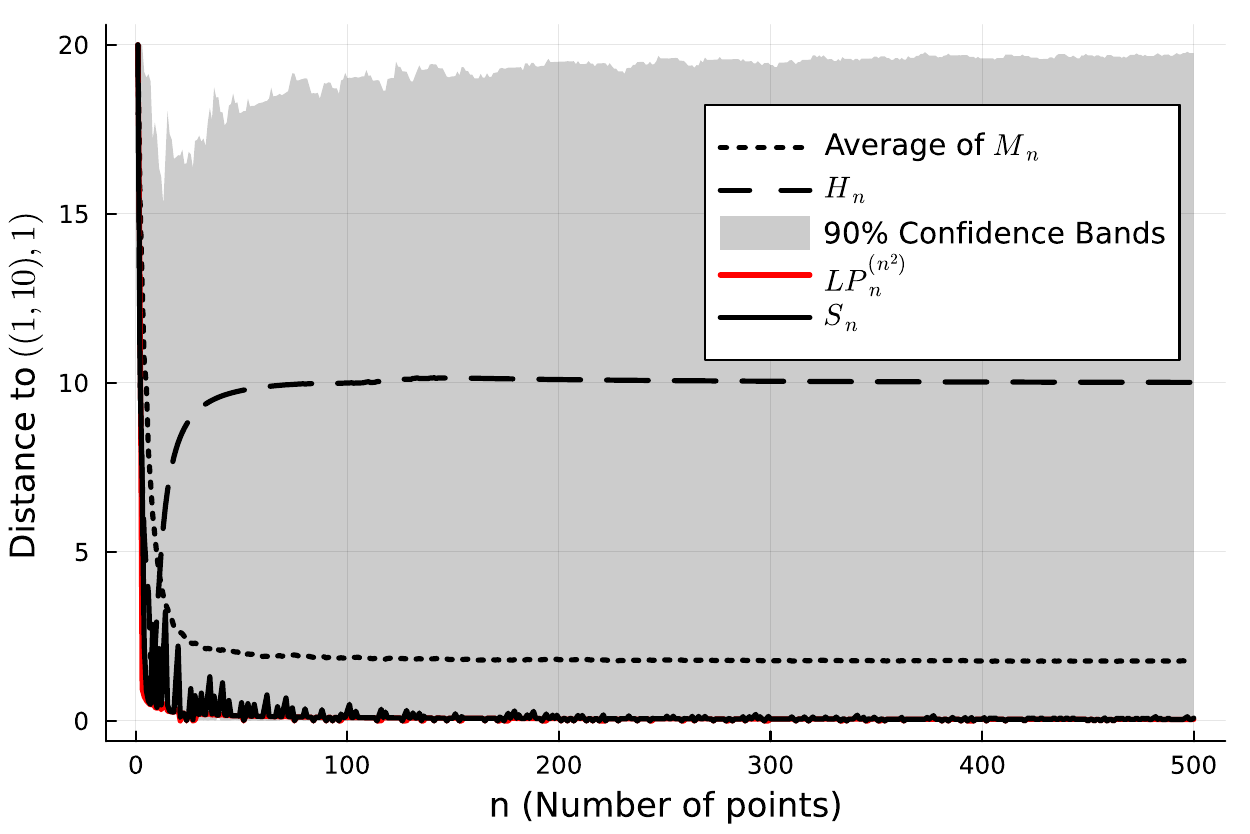}
  \caption{35\% Contamination.}
\end{subfigure}
\caption{Metric distance of $M_{n}, S_{n}$ and $H_{n}$ to the point $((1,10), 1)$ at various levels of contamination.}
\end{figure}

\begin{figure}
\centering
\includegraphics[width=0.3\textwidth]{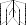}
\caption{Open book with 4 sheets of dimension 2. Each dashed line extends to infinity. To go from one point $p$ to a point $q$ in another sheet, one has to cross the spine.}
\label{fig:book-sketch}
\end{figure}

As a second example, we consider \textit{open books} (cf. \cite{hotz-sticky}). Intuitively speaking, open books are half-spaces glued together at a common \textit{spine}. Formally, let $d$ be an integer, and set
\begin{equation*}
  H_{+}^{d} = [0,\infty) \times \mathbb{R}^{d}.
\end{equation*}
$H_{+}^{d}$ can be thought of as a subset of $\mathbb{R}^{d+1}$ with boundary $\mathbb{R}^{d}$ (which is identified with $\{0\}\times \mathbb{R}^{d}$) and interior $(0,\infty) \times \mathbb{R}^{d}$. 
As such, $H_{+}^{d}$ inherits a subspace topology from $\mathbb{R}^{d+1}$. 
The open book $B_{k}^{d}$ with $k$ sheets of dimension $d+1$ is the disjoint union of $k$ copies of $H_{+}^{d}$ modulo identification of their boundaries, i.e., $B_{k}^{d} = (H_{+}\times \{1,\dots,k\})\,\big/\mkern-4.5mu\sim$, where $((t,x),j) \sim ((s,y),k)$ if and only if $t = s = 0$ and $x=y$. 

For the simulation, we look at a sequence of points $x_{n} = ((1+2/n, 10-1/\sqrt{n}),1)$ in $B_{3}^{1}$. 
This sequence clearly converges to $x = ((1,10),1)$. 
Again we use our contamination model.
The noise $y = ((1,10), s)$ has the same coordinates $(1,10)$, but may lie in a different, randomly selected, sheet $s=1,2,3$.

It is well known, that the classical stochastic limit theorems on open books defy Euclidean intuition.
Means in these spaces are \textit{sticky}, as Hotz et al. \cite{hotz-sticky} have noted, i.e., they are much more robust to outliers when compared to the classical Euclidean mean.

This has interesting consequences for the asymptotic recovery of means. 
Up to almost one-third of the data can be contaminated, and we can still recover the limit $x$. 
Both $S_{n}$ and $LP_{n}^{(n^{2})}$ are still able to recover $x$, even if more than one-third of the data is contaminated. 
In the case of positive definite matrices analyzed earlier, we saw, that $M_{n}$ and $S_{n}$ diverge, when a certain degree of contamination (about 5\%) is reached. 
Here $S_{n}$ and $LP_{n}^{(n^{2})}$ allow for twice as high a degree of contamination as $M_{n}$.
This means that ignoring information is much more detrimental in open books than it is in the case of symmetric, positive definite matrices.

Furthermore, the behavior of Hansen's mean $H_{n}$ differs significantly from the case of symmetric, positive definite matrices. When there is no contamination $H_{n}$ recovers $x$. However, once 1-2\% of the dataset is contaminated, we observe a phase transition. In this case Hansen's $H_{n}$ does not recover the true mean anymore. Judging from the simulations alone, it is unclear whether $H_{n}$ converges for low degrees of contamination. 
With high degrees of contamination, Hansen's $H_{n}$ seems to converge to the spine, i.e., to the point $((0,1),1)$. 
This is in line with the sticky behavior of means observed by Hotz et al. \cite{hotz-sticky}.   

\section{Discussion}\label{sec:disc}
In Theorem \ref{thm:slln-main} we have seen a (strong) law of large numbers for heteroscedastic data. 
Whether the Fr\'echet mean satisfies such a law is still open. 
The most recent convergence results for the Fr\'echet mean on Hadamard spaces \cite{le-gouic-fast-convergence, brunel, escande2023} either impose geometric assumptions on the underlying space, or the random variables themselves. 
To the best of our knowledge, there is no strong law of large numbers for the Fr\'echet mean of a general heteroskedastic sample of square-integrable, Hadamard space-valued random variables, as there is for the inductive mean. 
However, since the assumptions of Theorem \ref{thm:limit-bn} are weaker than the ones in Theorem \ref{thm:slln-main}, we know that, if such a law exists, the Fr\'echet mean must converge to the same limit as the inductive mean. 

In light of this, the recent results of Brunel and Serres \cite{brunel}, where it is shown that $S_n$ and $b_n$ are close with high probability, and the no-dice theorems of \cite{holbrook-no-dice} and \cite{lim-palfia}, one can not help but ask, what are the most general assumptions under which $\mathrm{E}(d(S_n,b_n)^2) \leq C/n$, for some $C>0$?
This question seems to be open in general.
Even if we restrict ourselves to the simpler case of a sequence $x_1,\dots,x_n$ of points in a Hadamard space, rather than sequences of Hadamard space-valued random variables, there does not seem to be a general estimate of $d(S_n,b_n)$.
Inspecting the proof of Lemma \ref{lemma:cesaro}, we observe that a condition like
\begin{equation*}
  \frac{1}{n}\sum_{k=1}^{n}d(x_{k},x)^{2} \to 0
\end{equation*}
is sufficient to guarantee that the inductive mean $S_{n}$ and the barycenter $b_{n}$ converge to the same point $x$.
This condition is of course weaker than the convergence of $x_{n} \to x$, but does in general not imply the weak convergence of $x_{n}$ to $x$ (see Section 3.1, \cite{dudley}).
Is weak convergence of the underlying sequence strong enough to guarantee that the inductive mean and the Fr\'echet mean agree asymptotically?

While the Fr\'echet mean is known to have a finite sample breakdown point of $1/n$, Propositions \ref{prop:stab-of-two-seq} -- \ref{prob:huber-L2}, and the results of Section \ref{sec:means} show that means in Hadamard spaces perform in general comparable to means in Euclidean spaces.
In light of the theoretical guarantees of our resampling scheme, and the simulation study, the inductive mean, although not symmetric, is quite robust with respect to data loss, data permutation, and noise. 
In the case of symmetric, positive definite matrices, Hansen's $H_{n}$ seems to converge against a point differing from the asymptotic \textit{center of mass}. 
This behavior is much more stable with respect to contaminated data than the convergence of the other procedures. Can Hansen's mean serve as a robust alternative to $S_{n}$ or $b_{n}$?

The simulations show, that the behavior of means crucially depends on the geometry of the underlying space.
Means of points close to the spine of open books are very stable, much more so than in the case of positive definite matrices.
One may observe that open books have a curvature of $-\infty$ along the spine. 
Curvature bounds are also crucial in recent laws of large numbers for the Fr\'echet mean \cite{brunel}. 
Can one characterize the robustness or \textit{stickiness} of means in terms of the local curvature of Hadamard spaces?

\section*{Acknowledgments}
Both authors are part of the Research Unit 5381 of the German Research Foundation. Georg Köstenberger is supported by the Austrian Science Fund (FWF): I 5485-N and Thomas Stark is supported by the Austrian Science Fund (FWF): I 5484-N. The authors would like to thank Moritz Jirak and Tatyana Krivobokova for organizing the research seminar at the Department for Statistics and Operations Research in Vienna and the guest speaker, Victor-Emmanuel Brunel, for a wonderful introduction to statistics on metric spaces and his feedback during the preparation of this manuscript.

\bibliographystyle{plain}
\bibliography{main}

\begin{thebibliography}{10}

\bibitem{entropy-ahidarcoutrix}
Adil Ahidar-Coutrix, Thibaut~Le Gouic, and Quentin Paris.
\newblock Convergence rates for empirical barycenters in metric spaces: curvature, convexity and extendable geodesics.
\newblock {\em Probability Theory and Related Fields}, 177:323--368, 2018.

\bibitem{alexandrov57}
Aleksandr~Danilovich Alexandrov.
\newblock {Über eine Verallgemeinerung der Riemannschen Geometrie}.
\newblock {\em Schriftenreihe des Forschinstituts fur Mathematik}, 1:33--84, 1957.

\bibitem{ando-li-mathias}
Tsuyoshi Ando, Chi-Kwong Li, and Roy Mathias.
\newblock Geometric means.
\newblock {\em Linear Algebra and its Applications}, 385:305--334, 2004.
\newblock Special Issue in honor of Peter Lancaster.

\bibitem{antezana2023}
Jorge Antezana, Eduardo Ghiglioni, and Demetrio Stojanoff.
\newblock Ergodic theorem in cat(0) spaces in terms of inductive means.
\newblock {\em Ergodic Theory and Dynamical Systems}, 43(5):1433–1454, 2023.

\bibitem{med1}
P.~Batchelor, M.~Moakher, D.~Atkinson, F.~Calamante, and A.~Connelly.
\newblock A rigorous framework for diffusion tensor calculus.
\newblock {\em Magnetic resonance in medicine : official journal of the Society of Magnetic Resonance in Medicine / Society of Magnetic Resonance in Medicine}, 53:221--5, 01 2005.

\bibitem{bacakbook}
Miroslav Ba\v{c}\'{a}k.
\newblock {\em Convex Analysis and Optimization in Hadamard Spaces}.
\newblock De Gruyter, Berlin, München, Boston, 2014.

\bibitem{bacak_challenges}
Miroslav Ba\v{c}\'{a}k.
\newblock Old and new challenges in {Hadamard} spaces.
\newblock {\em Jpn. J. Math.}, 18(2):115--168, 2023.

\bibitem{Bhathia-Holbroo}
Rajendra Bhatia and John Holbrook.
\newblock Riemannian geometry and matrix geometric means.
\newblock {\em Linear Algebra and its Applications}, 413(2):594--618, 2006.
\newblock Special Issue on the 11th Conference of the International Linear Algebra Society.

\bibitem{bhattacharya-patrangenaru-2003}
R.~Bhattacharya and V.~Patrangenaru.
\newblock {Large sample theory of intrinsic and extrinsic sample means on manifolds}.
\newblock {\em The Annals of Statistics}, 31(1):1 -- 29, 2003.

\bibitem{Billera-phylogenetic-trees}
Louis~J. Billera, Susan~P. Holmes, and Karen Vogtmann.
\newblock Geometry of the space of phylogenetic trees.
\newblock {\em Advances in Applied Mathematics}, 27(4):733--767, 2001.

\bibitem{bini2013}
Dario~A. Bini and Bruno Iannazzo.
\newblock Computing the karcher mean of symmetric positive definite matrices.
\newblock {\em Linear Algebra and its Applications}, 438(4):1700--1710, 2013.
\newblock 16th ILAS Conference Proceedings.

\bibitem{bridson-haefliner-book}
Martin~R. Bridson and Andr\'{e} H{\"a}fliger.
\newblock {\em Metric Spaces of Non-Positive Curvature}.
\newblock Grundlehren der mathematischen Wissenschaften. Springer Berlin Heidelberg, 2011.

\bibitem{brunel}
Victor-Emmanuel Brunel and Jordan Serres.
\newblock Concentration of empirical barycenters in metric spaces.
\newblock {\em Preprint, arXiv:12303.01144}, 2023.

\bibitem{burago-ferleger-kononenko}
Dmitry Burago, Serge Ferleger, and Alexey Kononenko.
\newblock Uniform estimates on the number of collisions in semi-dispersing billiards.
\newblock {\em Annals of Mathematics}, 147(3):695--708, 1998.

\bibitem{choi-toeplitz}
Byoung~Jin Choi and Un~Cig Ji.
\newblock Toeplitz lemma in geodesic metric space and convergence of inductive means.
\newblock {\em Journal of Mathematical Analysis and Applications}, 465(2):713--722, 2018.

\bibitem{dudley}
Richard~M. Dudley.
\newblock {\em Real Analysis and Probability}.
\newblock Cambridge Studies in Advanced Mathematics. Cambridge University Press, 2002.

\bibitem{efron79}
B.~Efron.
\newblock Bootstrap methods: Another look at the jackknife.
\newblock {\em The Annals of Statistics}, 7(1):1--26, 1979.

\bibitem{es-sahib-heinich}
Aziz Es-Sahib and Henri Heinich.
\newblock Barycentre canonique pour un espace m\'etrique \`a courbure n\'egative.
\newblock {\em S\'eminaire de probabilit\'es de Strasbourg}, 33:355--370, 1999.

\bibitem{escande2023}
Paul Escande.
\newblock On the concentration of the minimizers of empirical risks, 2023.

\bibitem{evans-jaffe-2024}
Steven~N. Evans and Adam~Q. Jaffe.
\newblock {Limit theorems for Fréchet mean sets}.
\newblock {\em Bernoulli}, 30(1):419 -- 447, 2024.

\bibitem{med2}
P.~Fletcher and S.~Joshi.
\newblock Riemannian geometry for the statistical analysis of diffusion tensor data.
\newblock {\em Signal Processing}, 87:250--262, 02 2007.

\bibitem{frechet-mean}
Maurice Fr\'echet.
\newblock Les \'el\'ements al\'eatoires de nature quelconque dans un espace distanci\'e.
\newblock {\em Annales de l'institut Henri Poincar\'e}, 10(4):215--310, 1948.

\bibitem{Funano}
Kei Funano.
\newblock {Rate of convergence of stochastic processes with values in $\mathbb{R}$-trees and Hadamard manifolds}.
\newblock {\em Osaka Journal of Mathematics}, 47(4):911 -- 920, 2010.

\bibitem{gauss1809theoria}
C.F. Gauss.
\newblock {\em Theoria motus corporum coelestium in sectionibus conicis solem ambientium}.
\newblock Carl Friedrich Gauss Werke. sumtibus Frid. Perthes et I. H. Besser, 1809.

\bibitem{gromov1987}
M.~Gromov.
\newblock {\em Hyperbolic Groups}, pages 75--263.
\newblock Springer New York, New York, NY, 1987.

\bibitem{gromov-random-walks}
M.~Gromov.
\newblock Random walk in random groups.
\newblock {\em Geom. Funct. Anal.}, 13(1):73--146, 2003.

\bibitem{Hansen}
Frank Hansen.
\newblock Regular operator mappings and multivariate geometric means.
\newblock {\em Linear Algebra and its Applications}, 461:123--138, 2014.

\bibitem{holbrook-no-dice}
John Holbrook.
\newblock No dice: A deterministic approach to the cartan centroid.
\newblock {\em Journal of the Ramanujan Mathematical Society}, 27, 12 2012.

\bibitem{hotz-sticky}
Thomas Hotz, Stephan Huckemann, Huiling Le, J.~S. Marron, Jonathan~C. Mattingly, Ezra Miller, James Nolen, Megan Owen, Vic Patrangenaru, and Sean Skwerer.
\newblock Sticky central limit theorems on open books.
\newblock {\em The Annals of Applied Probability}, 23(6):2238--2258, 2013.

\bibitem{huber64}
Peter~J. Huber.
\newblock Robust estimation of a location parameter.
\newblock {\em The Annals of Mathematical Statistics}, 35(1):73--101, 1964.

\bibitem{huber65}
Peter~J. Huber.
\newblock A robust version of the probability ratio test.
\newblock {\em The Annals of Mathematical Statistics}, 36(6):1753--1758, 1965.

\bibitem{huckemann}
S.~F. Huckemann.
\newblock {Intrinsic inference on the mean geodesic of planar shapes and tree discrimination by leaf growth}.
\newblock {\em The Annals of Statistics}, 39(2):1098 -- 1124, 2011.

\bibitem{kim}
Sejong Kim, Hosoo Lee, and Yongdo Lim.
\newblock A fixed point mean approximation to the cartan barycenter of positive definite matrices.
\newblock {\em Linear Algebra and its Applications}, 496:420--437, 2016.

\bibitem{lawson-lim}
Jimmie Lawson and Yongdo Lim.
\newblock Weighted means and karcher equations of positive operators.
\newblock {\em Proceedings of the National Academy of Sciences}, 110(39):15626--15632, 2013.

\bibitem{lawson-lim-monotonic}
Jimmie~D. Lawson and Yongdo Lim.
\newblock Monotonic properties of the least squares mean.
\newblock {\em Mathematische Annalen}, 351:267--279, 2010.

\bibitem{le-gouic-fast-convergence}
Thibaut Le~Gouic, Quentin Paris, Philippe Rigollet, and Austin~J. Stromme.
\newblock Fast convergence of empirical barycenters in {A}lexandrov spaces and the {W}asserstein space.
\newblock {\em J. Eur. Math. Soc. (JEMS)}, 25(6):2229--2250, 2023.

\bibitem{lim-palfia}
Yongdo Lim and Miklós Pálfia.
\newblock Weighted deterministic walks for the least squares mean on hadamard spaces.
\newblock {\em Bulletin of the London Mathematical Society}, 46(3):561--570, 2014.

\bibitem{loeve}
Michel Lo{\`e}ve.
\newblock {\em Probability Theory}.
\newblock Graduate texts in mathematics. Springer, 1977.

\bibitem{moakher2005}
Maher Moakher.
\newblock A differential geometric approach to the geometric mean of symmetric positive-definite matrices.
\newblock {\em SIAM J. Matrix Analysis Applications}, 26:735--747, 01 2005.

\bibitem{navas-ergodic}
Andr\'{e}s Navas.
\newblock An ${L}^1$ ergodic theorem with values in a non-positively curved space via a canonical barycenter map.
\newblock {\em Ergodic Theory and Dynamical Systems}, 33(2):609–623, 2013.

\bibitem{pennec2019riemannian}
Xavier Pennec, Stefan Sommer, and Tom Fletcher.
\newblock {\em Riemannian geometric statistics in medical image analysis}.
\newblock Academic Press, 2019.

\bibitem{qing-rafi-random-walks}
Yulan Qing and Kasra Rafi.
\newblock Sublinearly morse boundary i: Cat(0) spaces.
\newblock {\em Advances in Mathematics}, 404:108442, 2022.

\bibitem{schötz18-entropy}
Christof Sch{\"o}tz.
\newblock Convergence rates for the generalized fr{\'e}chet mean via the quadruple inequality.
\newblock {\em Electronic Journal of Statistics}, 2018.

\bibitem{schoetz}
Christof Schötz.
\newblock Strong laws of large numbers for generalizations of fréchet mean sets.
\newblock {\em Statistics}, 56(1):34--52, 2022.

\bibitem{sen-singer}
Pranab~K. Sen and Julio~M. Singer.
\newblock {\em Large Sample Methods in Statistics: An Introduction with Applications}.
\newblock Chapman \& Hall/CRC Texts in Statistical Science. Taylor \& Francis, 1994.

\bibitem{kiyoshi}
Kiyoshi Shiga.
\newblock Hadamard manifolds.
\newblock In {\em Geometry of geodesics and related topics ({T}okyo, 1982)}, volume~3 of {\em Adv. Stud. Pure Math.}, pages 239--281. North-Holland, Amsterdam, 1984.

\bibitem{sturm-martingale}
Karl-Theodor Sturm.
\newblock {Nonlinear martingale theory for processes with values in metric spaces of nonpositive curvature}.
\newblock {\em The Annals of Probability}, 30(3):1195 -- 1222, 2002.

\bibitem{sturm}
Karl-Theodor Sturm.
\newblock Probability measures on metric spaces of nonpositive curvature.
\newblock {\em Contemp. Math.}, 338, 01 2003.

\bibitem{slln-centroid-cpt-metric-space}
Harald Sverdrup-Thygeson.
\newblock {Strong Law of Large Numbers for Measures of Central Tendency and Dispersion of Random Variables in Compact Metric Spaces}.
\newblock {\em The Annals of Statistics}, 9(1):141 -- 145, 1981.

\bibitem{yokota}
Takumi Yokota.
\newblock {Convex functions and barycenter on CAT(1)-spaces of small radii}.
\newblock {\em Journal of the Mathematical Society of Japan}, 68(3):1297 -- 1323, 2016.

\bibitem{zhang2016}
Hongyi Zhang and Suvrit Sra.
\newblock First-order methods for geodesically convex optimization.
\newblock {\em Proceedings of Machine Learning Research}, 02 2016.

\bibitem{ziezold1977}
Herbert Ziezold.
\newblock {\em On Expected Figures and a Strong Law of Large Numbers for Random Elements in Quasi-Metric Spaces}, pages 591--602.
\newblock Springer Netherlands, Dordrecht, 1977.

\end{thebibliography}

\end{document}